\documentclass[1pt]{amsart}

\usepackage{amsmath,amsthm, amscd, amssymb, amsfonts}

\usepackage{hhline}

\newcommand{\Ind}{\operatorname{Ind}}
\newcommand{\ord}{\operatorname{ord}}
\newcommand{\oper}{\#}

\DeclareMathOperator*{\prode}{\boxtimes}

\newcommand\toba{{\mathfrak B }}
\newcommand\trasp{\pi}
\newcommand{\gr}{\operatorname{gr}}
\newcommand{\trid}{\triangleright}

\newcommand\compvert[2]{\genfrac{}{}{-1pt}{0}{#1}{#2}}

\newcommand\mvert[2]{\begin{tiny}\begin{matrix}#1\vspace{-4pt}\\#2\end{matrix}\end{tiny}}
\DeclareMathOperator*{\Tim}{\times}
\newcommand{\Times}[2]{\sideset{_#1}{_#2}\Tim}

\newcommand{\cx}{{\daleth}}

\newcommand{\lu}{{\nu}}
\newcommand{\ld}{{\ell}}

\newcommand{\abofe}{{A(\lu, \ld)}}

\newcommand{\bbofe}{{B(\lu, \ld)}}

\newcommand{\proy}{\Pi}

\newcommand{\Lc}{{\mathcal L}}
\newcommand{\Y}{{\mathcal Y}}
\newcommand{\R}{{\mathcal R}}
\newcommand{\W}{{\mathcal W}}
\newcommand{\Zc}{{\mathcal Z}}
\newcommand{\daga}{{\dagger}}

\newcommand{\acts}{{\rightharpoondown}}

\newcommand{\ku}{\mathbb C}
\newcommand{\K}{{\mathcal K}}
\newcommand{\Z}{{\mathbb Z}}
\newcommand{\N}{{\mathbb N}}
\newcommand{\I}{{\mathbb I}}

\newcommand{\G}{{\mathbb G}}
\newcommand{\M}{{\mathcal M}}
\newcommand{\Q}{{\mathcal Q}}
\newcommand{\F}{{\mathcal F}}
\newcommand{\C}{{\mathcal C}}
\newcommand{\D}{{\mathcal D}}
\newcommand{\Ec}{{\mathbf E}}
\newcommand{\Ee}{{\mathcal E}}
\newcommand{\Kc}{{\mathbf K}}

\newcommand{\uno}{{\bf 1}}
\newcommand{\uv}{{\bf 1_v}}
\newcommand{\uh}{{\bf 1_h}}
\newcommand{\m}{\mathcal{M}}
\newcommand{\n}{\mathcal{N}}

\newcommand{\Dc}{{\mathbf D}}
\newcommand{\B}{{\mathcal B}}
\newcommand{\T}{{\mathcal T}}
\newcommand{\Hc}{{\mathcal H}}
\newcommand{\Vc}{{\mathcal V}}
\newcommand{\Pc}{{\mathcal P}}
\newcommand{\Oc}{{\mathcal O}}
\newcommand{\ydh}{{}^H_H\mathcal{YD}}

\newcommand{\Ss}{{\mathcal S}}
\newcommand{\Vect}{\operatorname{Vec}}
\newcommand{\End}{\operatorname{End}}
\newcommand{\Aut}{\operatorname{Aut}}
\newcommand{\Int}{\operatorname{Int}}
\newcommand{\Ext}{\operatorname{Ext}}
\newcommand\tr{\operatorname{tr}}
\newcommand\Rep{\operatorname{Rep}}

\newcommand\card{\operatorname{card}}
\newcommand{\unosigma}{{\bf 1}^{\sigma}}
\newcommand\tot{\operatorname{tot}}
\newcommand\sgn{\operatorname{sgn}}
\newcommand\ad{\operatorname{ad}}
\newcommand\Hom{\operatorname{Hom}}
\newcommand\opext{\operatorname{Opext}}
\newcommand\Tot{\operatorname{Tot}}
\newcommand\Map{\operatorname{Map}}
\newcommand{\fde}{{\triangleright}}
\newcommand{\fiz}{{\triangleleft}}
\newcommand{\wC}{\widehat{C}}
\newcommand{\la}{\langle}
\newcommand{\ra}{\rangle}
\newcommand{\Comod}{\mbox{\rm Comod\,}}
\newcommand{\Mod}{\mbox{\rm Mod\,}}
\newcommand{\Sg}{{\mathfrak S}}
\newcommand{\funcion}{\mathfrak p}
\newcommand{\tauo}{{\widehat\tau}}
\newcommand{\prin}{t}
\newcommand{\fin}{b}
\newcommand{\pri}{r}
\newcommand{\fine}{l}
\newcommand\rh{\sim_{H}}
\newcommand\rv{\sim_{V}}
\newcommand\rd{\sim_{D}}

\newcommand\V{\operatorname{Vec}}

\theoremstyle{plain}

\renewcommand{\baselinestretch}{1.2}

\title{The character tables of  centralizers in  Sporadic Simple Groups of ${\rm McL}$}
\author{ \small Shouchuan Zhang,    \ \  Jieqiong He,\ \ Guichao Wu
}
\address{ Mathematics Department of Hunan University,\newline \indent Changsha China,
410082, E-mail: z9491@yahoo.com.cn }

\date{}

\begin{document}
\newtheorem{Proposition}{\quad Proposition}[section]
\newtheorem{Theorem}{\quad Theorem}
\newtheorem{Definition}[Proposition]{\quad Definition}
\newtheorem{Corollary}[Proposition]{\quad Corollary}
\newtheorem{Lemma}[Proposition]{\quad Lemma}
\newtheorem{Example}[Proposition]{\quad Example}

\maketitle \addtocounter{section}{-1}

\numberwithin{equation}{section}

\date{}

\begin {abstract}  To classify the finite dimensional pointed Hopf
algebras with   $G= {\rm McL}$ we obtain the representatives of
conjugacy classes of $G$ and all character tables of centralizers of
these representatives by means of software {\rm GAP}.


keywords: {\rm GAP}, Hopf algebra, sporadic simple group, character.
\end {abstract}

\section{Introduction}\label {s0}

This article is to contribute to the classification of
finite-dimensional complex pointed Hopf algebras  with sporadic
simple group $G= {\rm McL}$

 Many papers are about the classification of finite dimensional
pointed Hopf algebras, for example,  \cite{ AS02, AS00, AS05, He06,
AHS08, AG03, AFZ08, AZ07,Gr00,  Fa07, AF06, AF07, ZZC04, ZCZ08,
ZZWCY08}. In these research  ones need  the centralizers and
character tables of groups. In this paper we obtain   the
representatives of conjugacy classes of   sporadic simple group $
{\rm McL}$, as well as  all character tables of centralizers of
these representatives by means of software {\rm GAP}.

\section {${\rm McL}$ }

In this section $G$ denotes the sporadic simple group  ${\rm McL}$.
We use a  representation of $G$ given in \cite {Atlas}.

\subsection {Program}

gap $>$

a:= $\left(\begin{array}{cccccccccccccccccccccc}
1&0&0&0&0&0&0&0&0&0&0&0&0&0&0&0&0&0&0&0&0&0\\
0&-1&0&-1&1&0&0&0&0&0&0&0&-1&-1&-1&-1&0&0&0&0&1&0\\
0&0&0&0&0&0&1&0&0&0&0&0&0&0&0&0&0&0&0&0&0&0\\
1&-1&1&0&-1&0&0&0&0&0&-1&-1&0&-1&0&0&0&0&0&0&-1&0\\
0&0&0&0&0&0&0&-1&0&0&1&1&-1&-1&-1&0&1&0&1&0&0&1\\
-1&0&0&-1&0&0&0&0&0&0&0&0&-1&0&0&0&0&-1&0&0&0&0\\
0&0&1&0&0&0&0&0&0&0&0&0&0&0&0&0&0&0&0&0&0&0\\
1&0&1&0&-1&0&0&1&0&0&-1&-1&1&0&1&0&0&0&0&0&-1&0\\
-1&1&0&0&0&0&0&1&0&1&0&0&1&1&1&1&0&0&0&1&-1&0\\
0&0&-1&0&0&-1&1&-1&1&0&0&0&-1&0&0&-1&0&-1&0&-1&1&0\\
0&0&0&0&0&-1&1&0&0&0&0&-1&1&1&1&0&0&-1&0&0&0&0\\
1&1&0&1&-1&1&0&0&0&0&0&1&1&0&0&1&0&1&0&0&-1&0\\
-1&1&-1&0&0&0&0&0&0&0&1&1&0&0&0&0&0&0&0&0&0&0\\
0&0&0&0&0&0&0&0&0&0&0&0&0&1&0&0&0&0&0&0&0&0\\
0&0&0&0&0&0&0&0&0&0&0&0&0&0&0&1&0&0&0&0&0&0\\
0&0&0&0&0&0&0&0&0&0&0&0&0&0&1&0&0&0&0&0&0&0\\
0&0&1&-1&0&0&-1&1&0&0&-1&-1&0&0&0&1&0&0&-1&0&-1&-1\\
-1&0&0&0&1&-1&0&0&0&0&0&0&0&1&0&0&0&0&0&0&1&0\\
0&0&-1&0&0&0&1&0&0&0&0&0&0&0&1&-1&0&0&1&0&0&0\\
1&0&1&0&0&1&-1&1&0&0&-1&-1&1&0&0&1&-1&1&-1&1&-1&-1\\
0&0&0&-1&0&0&0&1&0&0&-1&-1&0&0&1&0&-1&0&-1&0&0&-1\\
-1&0&0&0&1&0&-1&0&0&0&1&1&-1&0&-1&0&0&0&0&0&1&1\\
\end{array}\right)$;;

 gap $>$

 b:=$\left(\begin{array}{cccccccccccccccccccccc}
   0& 0& 0& 0& 0& 0& 0& -1& 0& 0& 1& 1& -1& -1& -1& 0& 1& 0& 1& 0& 0& 1 \\
   0& 0& 0& 0& -1& 1& 0& -1& 0& 0& 0& 1& 0& 0& 0& 1& 0& 0& 0& 0& 0& 0 \\
   0& 0& 0& 0& 0& 1& 0& 0& 0& 0& 0& 0& 0& 0& 0& 0& 0& 0& 0& 0& 0& 0 \\
   0& 0& 0& 0& 0& 0& 0& 1& 0& 0& 0& -1& 1& 0& 1& 0& 0& 0& 0& 0& -1& 0 \\
   0& 0& 0& 0& 0& 0& 0& 0& 1& 0& 0& 0& 0& 0& 0& 0& 0& 0& 0& 0& 0& 0 \\
   0& -1& 0& -1& 0& 0& 0& 0& 0& 0& 0& -1& 0& 0& 0& 0& 0& -1& 0& 0& 0& 0 \\
   0& 0& 0& 0& 0& 0& 0& 0& 0& 0& 0& 0& 1& 0& 0& 0& 0& 0& 0& 0& 0& 0 \\
   -1& 0& -1& 0& 1& -1& 0& 0& 0& 0& 1& 0& 0& 0& 0& -1& 0& 0& 1& 0& 1& 1 \\
   0& 0& 0& 0& 0& 0& 0& 0& 0& 1& 0& 0& 0& 0& 0& 0& 0& 0& 0& 0& 0& 0 \\
   0& 0& 0& 0& 0& 1& -1& 0& 0& 0& 0& 0& 0& 0& 0& 1& 0& 0& 0& 1& -1& 0 \\
   0& 0& 0& 0& 1& 0& 0& 0& 0& 0& 0& 0& 0& 0& 0& 0& 0& 0& 0& 0& 0& 0 \\
   -1& 0& 0& 0& 1& -1& 0& 0& 0& 0& 1& 0& -1& 0& -1& -1& 1& 0& 1& 0& 1& 1 \\
   0& 0& 0& 0& 0& 0& 0& 0& 0& 0& 0& 0& 0& 1& 0& 0& 0& 0& 0& 0& 0& 0 \\
   0& 0& 0& 0& 0& 0& -1& 0& 0& 0& 0& 0& -1& -1& -1& 0& 0& 0& 0& 0& 0& 0 \\
   0& 0& 0& 0& 0& 0& 1& 0& 0& 0& 0& 0& 0& 0& 0& 0& 0& 0& 0& 0& 0& 0 \\
   1& 0& 0& 1& 0& 0& 0& 0& 0& -1& 0& 0& 1& 0& 0& -1& -1& 1& 0& 0& 0& 0 \\
   -1& 0& 0& 0& 0& -1& 0& 1& 0& 1& 0& 0& 0& 0& 0& 0& 0& 0& 0& 0& 0& 1 \\
   0& 0& 0& 0& 0& 0& 0& 0& 0& 0& -1& 0& 0& 0& 0& 0& -1& 0& -1& 0& 0& -1 \\
   0& 0& 0& 0& 0& 1& -1& 0& 0& 0& 0& 0& 0& 0& 0& 0& -1& 1& 0& 1& 0& 0 \\
   -1& 1& 0& 0& 0& 0& 0& 0& 0& 1& 1& 1& 0& 0& 0& 1& 1& 0& 0& 0& 0& 1 \\
   0& 0& 0& 0& 0& 0& 0& 0& 0& 0& 0& 0& 0& 0& 0& 0& 0& 0& 0& 0& 0& 1 \\
   1& 0& 0& 0& -1& 0& 1& -1& 0& -1& 0& 0& 0& 0& 0& 0& 0& -1& 0& -1& 0&
   -1\\
\end{array}\right)$;;

gap$>$ G:=Group(a,b);;

gap$>$ ccl:=ConjugacyClasses(G);;

gap$>$ q:=NrConjugacyClasses(G);;Display (q);

gap$>$ for i in [1..q] do

$>$ s:=Representative(ccl[i]);;Display(s);

$>$ od;

gap $>$  ccl:=ConjugacyClasses(G);;

gap $>$ q:=NrConjugacyClasses(G);; Display (q);

gap $>$  for i in [1..q] do

$>$ r:=Order(Representative(ccl[i])); Display(r);

$>$ od;

gap $>$  for i in [1..q] do

$>$ s:=Representative(ccl[i]);; cen:=Centralizer(G,s);;

$>$cl := ConjugacyClasses(cen);; t := NrConjugacyClasses(cen);;

$>$ for j in [1..t] do

$>$ if (s in cl[j]) then

$>$ Display(j);break;

$>$ fi;od;

$>$ od;

gap$>$ for i in [1..q] do

$>$ s:=Representative(ccl[i]);;cen:=Centralizer(G,s);;

$>$ cl:=ConjugacyClasses(cen);;t:=NrConjugacyClasses(cen);;

$>$ for j in [1..r] do

$>$ if (s$^j$ in ccl[i] and ((s$^j $= s) = false)) then

$>$ Print($" i=$",i,"$ AND  j=$",j, "$\n"$);

$>$ fi;od;od;

gap$>$ for i in [1..q] do

$>$ s:=Representative(ccl[i]);;

$>$ r:=Order(s);Display(r);

$>$ for k in [1..r] do

$>$ if (($(s^(k^2) = s^k)$ = false) and ($(s^(k^2) = s)$ = false))
then

$>$ Print($" i=$",i,"$ AND k =$ ",k,"$\n"$);

$>$ fi;od;od;

gap$>$ for i in [1..q] do

$>$ s:=Representative(ccl[i]);;cen:=Centralizer(G,s);;

$>$ cl:=ConjugacyClasses(cen);;t:=NrConjugacyClasses(cen);;

$>$ for j in [1..r] do

$>$ if (s$^j$ in ccl[i] and ($(s^j = s)$ = false)) and (($(s^(j^2) =
s^j)$ = false) and ($(s^(j^2) = s)$ = false)) then

$>$ Print($" i=$",i,"$ AND  j=$",j, "$\n"$);

$>$ fi;od;od;

gap$>$  s1:=Representative(ccl[1]);;cen1:=Centralizer(G,s1);;Display
(cen1);

gap$>$ cl1:=ConjugacyClasses(cen1);

gap$>$ char:=CharacterTable (G);;

gap$>$ Display (char);

gap$>$  s1:=Representative(ccl[2]);;cen1:=Centralizer(G,s1);;Display
(cen1);

gap$>$ cl1:=ConjugacyClasses(cen1);char:=CharacterTable (cen1);;

gap$>$ Display (char);

gap$>$  s1:=Representative(ccl[3]);;cen1:=Centralizer(G,s1);;Display
(cen1);

gap$>$ cl1:=ConjugacyClasses(cen1);char:=CharacterTable (cen1);;

gap$>$ Display (char);

gap$>$  s1:=Representative(ccl[4]);;cen1:=Centralizer(G,s1);;Display
(cen1);

gap$>$ cl1:=ConjugacyClasses(cen1);char:=CharacterTable (cen1);;

gap$>$ Display (char);

gap$>$  s1:=Representative(ccl[5]);;cen1:=Centralizer(G,s1);;Display
(cen1);

gap$>$ cl1:=ConjugacyClasses(cen1);char:=CharacterTable (cen1);;

gap$>$ Display (char);

gap$>$  s1:=Representative(ccl[6]);;cen1:=Centralizer(G,s1);;Display
(cen1);

gap$>$ cl1:=ConjugacyClasses(cen1);char:=CharacterTable (cen1);;

gap$>$ Display (char);

gap$>$  s1:=Representative(ccl[7]);;cen1:=Centralizer(G,s1);;Display
(cen1);

gap$>$ cl1:=ConjugacyClasses(cen1);char:=CharacterTable (cen1);;

gap$>$ Display (char);

gap$>$  s1:=Representative(ccl[8]);;cen1:=Centralizer(G,s1);;Display
(cen1);

gap$>$ cl1:=ConjugacyClasses(cen1);char:=CharacterTable (cen1);;

gap$>$ Display (char);

gap$>$  s1:=Representative(ccl[9]);;cen1:=Centralizer(G,s1);;Display
(cen1);

gap$>$ cl1:=ConjugacyClasses(cen1);char:=CharacterTable (cen1);;

gap$>$ Display (char);

gap$>$
s1:=Representative(ccl[10]);;cen1:=Centralizer(G,s1);;Display
(cen1);

gap$>$ cl1:=ConjugacyClasses(cen1);char:=CharacterTable (cen1);;

gap$>$ Display (char);

gap$>$
s1:=Representative(ccl[11]);;cen1:=Centralizer(G,s1);;Display
(cen1);

gap$>$ cl1:=ConjugacyClasses(cen1);char:=CharacterTable (cen1);;

gap$>$ Display (char);

gap$>$
s1:=Representative(ccl[12]);;cen1:=Centralizer(G,s1);;Display
(cen1);

gap$>$ cl1:=ConjugacyClasses(cen1);char:=CharacterTable (cen1);;

gap$>$ Display (char);

gap$>$
s1:=Representative(ccl[13]);;cen1:=Centralizer(G,s1);;Display
(cen1);

gap$>$ cl1:=ConjugacyClasses(cen1);char:=CharacterTable (cen1);;

gap$>$ Display (char);

gap$>$
s1:=Representative(ccl[14]);;cen1:=Centralizer(G,s1);;Display
(cen1);

gap$>$ cl1:=ConjugacyClasses(cen1);char:=CharacterTable (cen1);;

gap$>$ Display (char);

gap$>$
s1:=Representative(ccl[15]);;cen1:=Centralizer(G,s1);;Display
(cen1);

gap$>$ cl1:=ConjugacyClasses(cen1);char:=CharacterTable (cen1);;

gap$>$ Display (char);

gap$>$
s1:=Representative(ccl[16]);;cen1:=Centralizer(G,s1);;Display
(cen1);

gap$>$ cl1:=ConjugacyClasses(cen1);char:=CharacterTable (cen1);;

gap$>$ Display (char);

gap$>$
s1:=Representative(ccl[17]);;cen1:=Centralizer(G,s1);;Display
(cen1);

gap$>$ cl1:=ConjugacyClasses(cen1);char:=CharacterTable (cen1);;

gap$>$ Display (char);

gap$>$
s1:=Representative(ccl[18]);;cen1:=Centralizer(G,s1);;Display
(cen1);

gap$>$ cl1:=ConjugacyClasses(cen1);char:=CharacterTable (cen1);;

gap$>$ Display (char);

gap$>$
s1:=Representative(ccl[19]);;cen1:=Centralizer(G,s1);;Display
(cen1);

gap$>$ cl1:=ConjugacyClasses(cen1);char:=CharacterTable (cen1);;

gap$>$ Display (char);

gap$>$
s1:=Representative(ccl[20]);;cen1:=Centralizer(G,s1);;Display
(cen1);

gap$>$ cl1:=ConjugacyClasses(cen1);char:=CharacterTable (cen1);;

gap$>$ Display (char);

gap$>$
s1:=Representative(ccl[21]);;cen1:=Centralizer(G,s1);;Display
(cen1);

gap$>$ cl1:=ConjugacyClasses(cen1);char:=CharacterTable (cen1);;

gap$>$ Display (char);

gap$>$
s1:=Representative(ccl[22]);;cen1:=Centralizer(G,s1);;Display
(cen1);

gap$>$ cl1:=ConjugacyClasses(cen1);char:=CharacterTable (cen1);;

gap$>$ Display (char);

gap$>$
s1:=Representative(ccl[23]);;cen1:=Centralizer(G,s1);;Display
(cen1);

gap$>$ cl1:=ConjugacyClasses(cen1);char:=CharacterTable (cen1);;

gap$>$ Display (char);

gap$>$
s1:=Representative(ccl[24]);;cen1:=Centralizer(G,s1);;Display
(cen1);

gap$>$ cl1:=ConjugacyClasses(cen1);char:=CharacterTable (cen1);;

gap$>$ Display (char);

\subsection {The character tables }

The order of $G$ is 898128000.

The generators of $G$ are:\\

$\left(
\begin{array}{cccccccccccccccccccccc}
1&0&0&0&0&0&0&0&0&0&0&0&0&0&0&0&0&0&0&0&0&0\\
0&-1&0&-1&1&0&0&0&0&0&0&0&-1&-1&-1&-1&0&0&0&0&1&0\\
0&0&0&0&0&0&1&0&0&0&0&0&0&0&0&0&0&0&0&0&0&0\\
1&-1&1&0&-1&0&0&0&0&0&-1&-1&0&-1&0&0&0&0&0&0&-1&0\\
0&0&0&0&0&0&0&-1&0&0&1&1&-1&-1&-1&0&1&0&1&0&0&1\\
-1&0&0&-1&0&0&0&0&0&0&0&0&-1&0&0&0&0&-1&0&0&0&0\\
0&0&1&0&0&0&0&0&0&0&0&0&0&0&0&0&0&0&0&0&0&0\\
1&0&1&0&-1&0&0&1&0&0&-1&-1&1&0&1&0&0&0&0&0&-1&0\\
-1&1&0&0&0&0&0&1&0&1&0&0&1&1&1&1&0&0&0&1&-1&0\\
0&0&-1&0&0&-1&1&-1&1&0&0&0&-1&0&0&-1&0&-1&0&-1&1&0\\
0&0&0&0&0&-1&1&0&0&0&0&-1&1&1&1&0&0&-1&0&0&0&0\\
1&1&0&1&-1&1&0&0&0&0&0&1&1&0&0&1&0&1&0&0&-1&0\\
-1&1&-1&0&0&0&0&0&0&0&1&1&0&0&0&0&0&0&0&0&0&0\\
0&0&0&0&0&0&0&0&0&0&0&0&0&1&0&0&0&0&0&0&0&0\\
0&0&0&0&0&0&0&0&0&0&0&0&0&0&0&1&0&0&0&0&0&0\\
0&0&0&0&0&0&0&0&0&0&0&0&0&0&1&0&0&0&0&0&0&0\\
0&0&1&-1&0&0&-1&1&0&0&-1&-1&0&0&0&1&0&0&-1&0&-1&-1\\
-1&0&0&0&1&-1&0&0&0&0&0&0&0&1&0&0&0&0&0&0&1&0\\
0&0&-1&0&0&0&1&0&0&0&0&0&0&0&1&-1&0&0&1&0&0&0\\
1&0&1&0&0&1&-1&1&0&0&-1&-1&1&0&0&1&-1&1&-1&1&-1&-1\\
0&0&0&-1&0&0&0&1&0&0&-1&-1&0&0&1&0&-1&0&-1&0&0&-1\\
-1&0&0&0&1&0&-1&0&0&0&1&1&-1&0&-1&0&0&0&0&0&1&1\\
\end{array}\right)$,

$\left(\begin{array}{cccccccccccccccccccccc}
   0& 0& 0& 0& 0& 0& 0& -1& 0& 0& 1& 1& -1& -1& -1& 0& 1& 0& 1& 0& 0& 1 \\
   0& 0& 0& 0& -1& 1& 0& -1& 0& 0& 0& 1& 0& 0& 0& 1& 0& 0& 0& 0& 0& 0 \\
   0& 0& 0& 0& 0& 1& 0& 0& 0& 0& 0& 0& 0& 0& 0& 0& 0& 0& 0& 0& 0& 0 \\
   0& 0& 0& 0& 0& 0& 0& 1& 0& 0& 0& -1& 1& 0& 1& 0& 0& 0& 0& 0& -1& 0 \\
   0& 0& 0& 0& 0& 0& 0& 0& 1& 0& 0& 0& 0& 0& 0& 0& 0& 0& 0& 0& 0& 0 \\
   0& -1& 0& -1& 0& 0& 0& 0& 0& 0& 0& -1& 0& 0& 0& 0& 0& -1& 0& 0& 0& 0 \\
   0& 0& 0& 0& 0& 0& 0& 0& 0& 0& 0& 0& 1& 0& 0& 0& 0& 0& 0& 0& 0& 0 \\
   -1& 0& -1& 0& 1& -1& 0& 0& 0& 0& 1& 0& 0& 0& 0& -1& 0& 0& 1& 0& 1& 1 \\
   0& 0& 0& 0& 0& 0& 0& 0& 0& 1& 0& 0& 0& 0& 0& 0& 0& 0& 0& 0& 0& 0 \\
   0& 0& 0& 0& 0& 1& -1& 0& 0& 0& 0& 0& 0& 0& 0& 1& 0& 0& 0& 1& -1& 0 \\
   0& 0& 0& 0& 1& 0& 0& 0& 0& 0& 0& 0& 0& 0& 0& 0& 0& 0& 0& 0& 0& 0 \\
   -1& 0& 0& 0& 1& -1& 0& 0& 0& 0& 1& 0& -1& 0& -1& -1& 1& 0& 1& 0& 1& 1 \\
   0& 0& 0& 0& 0& 0& 0& 0& 0& 0& 0& 0& 0& 1& 0& 0& 0& 0& 0& 0& 0& 0 \\
   0& 0& 0& 0& 0& 0& -1& 0& 0& 0& 0& 0& -1& -1& -1& 0& 0& 0& 0& 0& 0& 0 \\
   0& 0& 0& 0& 0& 0& 1& 0& 0& 0& 0& 0& 0& 0& 0& 0& 0& 0& 0& 0& 0& 0 \\
   1& 0& 0& 1& 0& 0& 0& 0& 0& -1& 0& 0& 1& 0& 0& -1& -1& 1& 0& 0& 0& 0 \\
   -1& 0& 0& 0& 0& -1& 0& 1& 0& 1& 0& 0& 0& 0& 0& 0& 0& 0& 0& 0& 0& 1 \\
   0& 0& 0& 0& 0& 0& 0& 0& 0& 0& -1& 0& 0& 0& 0& 0& -1& 0& -1& 0& 0& -1 \\
   0& 0& 0& 0& 0& 1& -1& 0& 0& 0& 0& 0& 0& 0& 0& 0& -1& 1& 0& 1& 0& 0 \\
   -1& 1& 0& 0& 0& 0& 0& 0& 0& 1& 1& 1& 0& 0& 0& 1& 1& 0& 0& 0& 0& 1 \\
   0& 0& 0& 0& 0& 0& 0& 0& 0& 0& 0& 0& 0& 0& 0& 0& 0& 0& 0& 0& 0& 1 \\
   1& 0& 0& 0& -1& 0& 1& -1& 0& -1& 0& 0& 0& 0& 0& 0& 0& -1& 0& -1& 0&
   -1\\
\end{array}\right)$.

The representatives of conjugacy classes of   $G$ are:\\
$s_1$=$\left(\begin{array}{cccccccccccccccccccccc}    1&  0&  0& 0&
0&  0&  0&  0&  0&  0&  0&  0&  0&  0& 0& 0&  0&  0&  0&  0&
 0&  0 \\
    0&  1&  0&  0&  0&  0&  0&  0&  0&  0&  0&  0&  0&  0&  0&  0&  0&  0&  0&  0&
 0&  0 \\
    0&  0&  1&  0&  0&  0&  0&  0&  0&  0&  0&  0&  0&  0&  0&  0&  0&  0&  0&  0&
 0&  0 \\
    0&  0&  0&  1&  0&  0&  0&  0&  0&  0&  0&  0&  0&  0&  0&  0&  0&  0&  0&  0&
 0&  0 \\
    0&  0&  0&  0&  1&  0&  0&  0&  0&  0&  0&  0&  0&  0&  0&  0&  0&  0&  0&  0&
 0&  0 \\
    0&  0&  0&  0&  0&  1&  0&  0&  0&  0&  0&  0&  0&  0&  0&  0&  0&  0&  0&  0&
 0&  0 \\
    0&  0&  0&  0&  0&  0&  1&  0&  0&  0&  0&  0&  0&  0&  0&  0&  0&  0&  0&  0&
 0&  0 \\
    0&  0&  0&  0&  0&  0&  0&  1&  0&  0&  0&  0&  0&  0&  0&  0&  0&  0&  0&  0&
 0&  0 \\
\end{array}\right)$,

$s_2$=$\left(\begin{array}{cccccccccccccccccccccc}   0&   0& 0&0& 0&
0&   0&  -1&   0&  -1&   0&  -1&   0&   1&   0&   0&
0&  -1&   0&   0&   0&  -1 \\
0&   0&   0&   1&   0&  -1&   0&   0&   0&   0&   0&   0&   0&   0&
0&  -1&0&   0&   1&   0&   0&   1 \\
0&   1&   1&   0&   0&   0&  -1&   1&   0&   0&   0&   0&   0&   0&
0&   1& 0&   0&  -1&   0&  -1&   0 \\
0&   0&   1&  -1&   1&  -1&   0&   1&   0&   0&   0&  -1&   0&   0&
0&   0&
  1&   0&   0&   0&   0&   0 \\
     0&   0&   0&   0&   0&   0&   1&  -1&  -1&  -1&   1&   0&   0&   1&   0&   0&
  1&  -1&   1&   0&   0&   0 \\
     0&   0&   0&   1&   0&   0&   0&   0&  -1&   0&   1&   1&   1&   0&  -1&   0&
  0&   1&   1&   1&   0&   1 \\
     0&   0&   0&   0&  -1&   0&   1&   0&  -1&   0&   0&   0&   1&   0&   1&   0&
  0&   0&   1&   0&  -1&   0 \\
    -1&   0&   0&  -1&   0&   0&   0&   0&   0&   1&   0&   0&   0&   0&   0&   1&
  0&   0&   0&   1&   0&   0 \\
\end{array}\right)$,

$s_3$=$\left(\begin{array}{cccccccccccccccccccccc}   1&   1& 0& 1&
-1&   1&   0&   0&   0&   0&   0&   1&   1&   0&   0&   1&
  0&   1&   0&   0&  -1&   0 \\
     1&  -1&   0&   1&   0&   0&   0&   0&   0&   0&   0&   0&   0&  -1&  -1&  -1&
  0&   1&   1&   0&   0&   1 \\
    -1&   0&   0&   0&   1&  -1&   0&   0&   0&   0&   1&   0&  -1&   0&  -1&  -1&
  1&   0&   1&   0&   1&   1 \\
     0&   0&   1&  -1&   0&  -1&   0&   0&   0&   0&   0&  -1&  -1&   0&   0&   0&
  1&  -1&   0&  -1&   0&   0 \\
     0&   0&   0&   0&   0&   0&   0&   0&  -1&  -1&   1&  -1&   1&   1&   0&   0&
  0&   0&   1&   1&  -1&   0 \\
    -1&   0&  -1&   0&   1&  -1&   0&   0&   0&   0&   1&   0&   0&   0&   0&  -1&
  0&   0&   1&   0&   1&   1 \\
     0&  -1&   1&  -1&   1&  -1&  -1&   1&   0&   0&   0&  -1&  -1&  -1&  -1&  -1&
  1&   0&   0&   0&   0&   1 \\
     0&   1&   0&   0&   0&  -1&   0&   1&   1&   1&  -1&   0&   0&   0&   1&   0&
  0&   0&  -1&  -1&   0&   0 \\
\end{array}\right)$,

$s_4$=$\left(\begin{array}{cccccccccccccccccccccc}   0&   1& 0& 1&
0&   0&  -1&   1&   1&   1&   0&   1&   0&  -1&   0&   0&
0&   1&   0&   0&   0&   1 \\
     1&   0&   1&   0&   0&   0&   0&   0&   0&   0&   0&  -1&   0&   0&   0&   0&
  1&   0&   0&   0&  -1&   0 \\
    -1&   0&   0&  -1&   1&  -1&   0&   0&   0&   0&   0&  -1&  -1&   0&   0&   0&
  1&  -1&   0&   0&   0&   0 \\
    -1&   0&   0&  -1&   0&   0&   0&   0&   0&   0&   0&   0&  -1&   0&   0&   0&
  0&  -1&   0&   0&   0&   0 \\
     0&   1&   0&   0&   0&   0&   0&   0&   0&   0&   0&   0&   0&   0&   0&   0&
  0&   0&   0&   0&   0&   0 \\
     0&   0&   1&  -1&   0&  -1&   0&   0&   0&   0&   0&  -1&  -1&   0&   0&   0&
  1&  -1&   0&  -1&   0&   0 \\
     0&   0&   0&   0&   0&   0&   0&   0&   0&   0&   0&   0&   0&   0&   0&   0&
  0&   0&   0&   0&   0&   1 \\
     0&   1&  -1&   1&  -1&   0&   0&   0&   1&   0&   0&   1&   0&   0&   1&   0&
 -1&   0&   0&  -1&   0&   0 \\
\end{array}\right)$,

$s_5$=$\left(\begin{array}{cccccccccccccccccccccc}     0&   0& 0& 0&
0&   1&  -1&   1&  -1&   0&   0&   0&   1&   0&   0&   0&
-1&   1&   0&   1&  -1&   0 \\
    0&   0&   0&   0&   0&   0&   0&   0&  -1&  -1&   1&   0&   1&   1&   0&   0&
  0&   0&   0&   0&   0&   0 \\
    0&   0&   0&   1&   0&  -1&   0&   0&   1&   0&   0&   0&   0&   0&   0&  -1&
  0&   0&   0&  -1&   1&   1 \\
    1&   0&   0&   1&  -1&   1&   0&  -1&   0&   0&   0&   1&   0&  -1&  -1&   0&
  0&   1&   0&   0&   0&   0 \\
    1&  -1&   0&   0&  -1&   0&   1&   0&   0&   0&  -1&   0&   0&  -1&   0&  -1&
  0&   0&   0&  -1&   0&   0 \\
    1&   0&  -1&   1&  -1&   1&   1&  -1&   0&  -1&   0&   0&   1&   1&   1&   0&
 -1&   0&   0&   0&   0&  -1 \\
    0&   1&   0&   1&  -1&   0&   0&   1&   1&   1&  -1&   1&   1&   0&   1&   1&
 -1&   1&  -1&   0&   0&   0 \\
    1&   0&   1&   0&   0&   1&  -1&   1&   0&   0&  -1&  -1&   1&   0&   0&   1&
 -1&   1&  -1&   1&  -1&  -1 \\
\end{array}\right)$,

$s_6$=$\left(\begin{array}{cccccccccccccccccccccc}   -1&   0& -1& 0&
0&   0&   1&  -1&   0&   0&   1&   0&   0&   1&   1&   0&
0&  -1&   1&   0&   0&   0 \\
   0&   0&   0&   0&  -1&   0&   0&   0&   0&   1&   0&   0&   0&  -1&   0&   0&
  1&   0&   1&   0&  -1&   1 \\
   0&   0&   1&   0&   1&  -1&  -1&   1&   0&   0&   0&  -1&   0&   0&  -1&   0&
  1&   0&   0&   0&   0&   1 \\
   0&   0&   0&   1&   0&   0&   0&   0&  -1&   0&   1&   1&   1&   0&  -1&   0&
  0&   1&   1&   1&   0&   1 \\
   0&  -1&   1&  -1&   0&   0&   0&   1&  -1&   1&   0&  -1&   0&  -1&   0&   0&
  1&   0&   1&   1&  -1&   1 \\
   0&   0&   0&   0&  -1&   1&   0&   0&  -1&   0&   0&   0&   1&   0&   0&   1&
  0&   0&   0&   1&  -1&   0 \\
   0&   0&   1&   0&   0&   0&  -1&   1&  -1&   0&   0&  -1&   1&   0&   0&   0&
  0&   1&   0&   1&  -1&   0 \\
  -1&   1&  -1&   1&   0&   0&   0&  -1&   0&  -1&   1&   1&   0&   1&   0&   0&
 -1&   0&   0&   0&   1&   0 \\
\end{array}\right)$,

$s_7$=$\left(\begin{array}{cccccccccccccccccccccc}    0&   0& 1& 0&
0&  -1&   0&   0&   0&   0&   0&   0&   0&   0&  -1&   0&
  1&   0&   0&   0&   0&   1 \\
   -1&   0&   0&   0&   0&  -1&   0&   1&   0&   1&   0&   0&   0&   0&   0&   0&
  0&   0&   0&   0&   0&   1 \\
    0&   0&   0&   0&   0&   0&   0&   1&   0&   0&   0&  -1&   1&   0&   1&   0&
  0&   0&   0&   0&  -1&   0 \\
    1&  -1&   1&  -1&   0&   1&   0&   0&  -1&  -1&   0&  -1&   0&   0&   0&   0&
  0&   0&   0&   0&  -1&  -1 \\
    1&  -1&   1&   0&  -1&   0&   1&   0&  -1&   0&   0&  -1&   1&   0&   0&   0&
  1&   0&   1&   0&  -1&   0 \\
    1&  -1&   1&   0&   0&  -1&   0&   1&   0&   0&  -1&  -1&   0&  -1&   0&  -1&
  0&   0&   0&  -1&   0&   0 \\
    0&  -1&   0&  -1&   0&   0&   0&   0&   0&   0&   0&  -1&   0&   0&   0&   0&
  0&  -1&   0&   0&   0&   0 \\
    0&  -1&   1&  -1&   1&   0&  -1&   0&  -1&  -1&   0&  -1&   0&   0&  -1&   0&
  0&   0&   0&   1&   0&   0 \\
\end{array}\right)$,

$s_8$=$\left(\begin{array}{cccccccccccccccccccccc}  0&  -1&   1& -1&
1&  -1&   0&   0&   0&   0&   0&  -1&  -1&   0&  -1&   0&
1&  -1&   0&   0&   0&   0 \\
    1&  -1&   1&   0&   0&   0&   0&   0&  -1&  -1&   0&  -1&   0&  -1&  -1&  -1&
  1&   0&   1&   0&  -1&   0 \\
    0&   0&  -1&   0&   0&   0&   1&   0&   0&   0&   0&   0&   0&   0&   1&  -1&
  0&   0&   1&   0&   0&   0 \\
    0&   1&  -1&   1&  -1&   1&   0&   0&   0&   1&   0&   1&   1&   0&   1&   1&
 -1&   1&   0&   1&  -1&   0 \\
    0&   0&   0&   1&   0&  -1&   0&   0&   0&   0&   0&   0&   0&   0&   0&  -1&
  0&   0&   1&   0&   0&   1 \\
    0&   1&   0&   1&   0&   0&  -1&   1&   1&   0&  -1&   0&   0&   0&   0&   0&
 -1&   1&  -1&   0&   0&   0 \\
    0&   0&  -1&   1&   0&   0&   1&  -1&  -1&  -1&   1&   0&   1&   1&   0&  -1&
  0&   0&   1&   0&   0&   0 \\
    0&   0&   0&   0&   0&   0&   0&   0&   1&   1&  -1&   0&   0&   0&   0&   1&
  0&   0&  -1&   0&   0&   0 \\
\end{array}\right)$,

$s_9$=$\left(\begin{array}{cccccccccccccccccccccc}  0&  -1&   0& 0&
0&   0&   0&  -1&  -1&   0&   1&   0&   0&   0&  -1&   0&
  1&   0&   1&   1&   0&   1 \\
   -1&   0&   0&  -1&   1&   0&  -1&   0&   0&   0&   0&   0&  -1&   0&  -1&   0&
  0&   0&  -1&   0&   1&   0 \\
    0&   0&   1&   0&   1&   0&  -1&   1&   0&   0&   0&   0&   0&   0&  -1&   0&
  0&   1&  -1&   0&   0&   0 \\
    1&  -1&   1&   0&   0&  -1&   0&   1&   0&   0&  -1&  -1&   0&  -1&   0&  -1&
  0&   0&   0&  -1&   0&   0 \\
    0&  -1&   0&   0&   0&   1&   0&   0&  -1&   0&   0&   0&   1&   0&   0&   0&
 -1&   1&   0&   1&   0&   0 \\
    0&   0&   0&   0&   0&   0&   0&   0&   0&   0&   0&   0&   1&   0&   0&   0&
  0&   0&   0&   0&   0&   0 \\
   -1&   0&   0&  -1&   1&  -1&   0&   0&   1&   0&   0&   0&  -1&   0&   0&   0&
  0&  -1&  -1&  -1&   1&   0 \\
   -1&   0&   0&   0&   0&  -1&   0&   0&   0&   0&   1&   0&   0&   0&   0&   0&
  1&  -1&   1&   0&   0&   1 \\
\end{array}\right)$,

$s_{10}$=$\left(\begin{array}{cccccccccccccccccccccc} 0& 0& 0& 0& 0&
-1&   0&   1&   1&   1&  -1&   0&   0&  -1&   0&   0&
0&   0&   0&   0&   0&   1 \\
     0&   0&   0&   0&   0&   0&   0&   0&  -1&  -1&   1&  -1&   1&   1&   0&   0&
  0&   0&   1&   1&  -1&   0 \\
     1&   0&   0&   1&  -1&   1&   0&  -1&  -1&  -1&   1&   0&   1&   0&   0&   0&
  0&   0&   1&   0&  -1&   0 \\
     0&   0&   0&   0&   0&   1&   0&   0&   0&   0&   0&   0&   0&   0&   0&   0&
  0&   0&   0&   0&   0&   0 \\
     1&  -1&   1&   0&   0&   0&   0&   0&  -1&  -1&   0&  -1&   0&  -1&  -1&  -1&
  1&   0&   1&   0&  -1&   0 \\
     1&  -1&   1&  -1&   0&   1&   0&   0&  -1&  -1&   0&  -1&   0&   0&   0&   0&
  0&   0&   0&   0&  -1&  -1 \\
     1&  -1&   1&   0&   0&   0&   0&   0&  -1&  -1&   0&  -1&   0&   0&  -1&  -1&
  0&   0&   0&   0&   0&   0 \\
     0&   0&   1&  -1&   0&   0&   0&   1&   0&   1&  -1&   0&   0&   0&   0&   1&
  0&   0&  -1&   0&   0&   0 \\
\end{array}\right)$,

$s_{11}$=$\left(\begin{array}{cccccccccccccccccccccc}  0&   0& 0& 0&
-1&   1&   0&   0&  -1&   0&   0&   0&   1&   0&   0&   1&
0&   0&   0&   1&  -1&   0 \\
    0&   1&   0&   1&  -1&   0&   0&   0&   1&   1&   0&   1&   0&   0&   0&   1&
  0&   0&   0&   0&   0&   1 \\
   -1&   0&   0&   0&   1&  -1&   0&   0&   0&   0&   1&   0&  -1&   0&  -1&  -1&
  1&   0&   1&   0&   1&   1 \\
    0&  -1&   0&   0&   1&   0&   0&  -1&   0&  -1&   0&   0&  -1&   0&  -1&  -1&
  0&   0&   0&   0&   1&   0 \\
    0&  -1&   1&   0&   0&   0&  -1&   1&   0&   1&  -1&   0&   0&  -1&  -1&   0&
  0&   1&   0&   1&   0&   1 \\
    1&   0&   1&   1&   0&   0&  -1&   0&   0&  -1&   0&   0&   0&   0&  -1&   0&
  0&   1&   0&   0&   0&   0 \\
    0&   0&   1&   0&   0&   0&   0&   0&   0&   0&   0&   0&   0&   0&   0&   0&
  0&   0&   0&   0&   0&   0 \\
    1&   0&   1&   0&   0&   0&   0&   0&   0&  -1&  -1&  -1&   0&   0&   0&   0&
  0&   0&  -1&  -1&   0&  -1 \\
\end{array}\right)$,

$s_{12}$=$\left(\begin{array}{cccccccccccccccccccccc}-1&   0& -1& 0&
0&   1&   0&  -1&  -1&   0&   1&   1&   0&   0&   0&   0&
0&   0&   1&   1&   0&   0 \\
     0&   0&   0&   1&   0&   0&   0&   0&   0&   0&   0&   0&   0&   0&   0&   0&
  0&   0&   0&   0&   0&   0 \\
    -1&   0&   1&  -1&   1&   0&  -1&   1&   0&   0&   0&  -1&   0&   1&   0&   1&
  0&   0&  -1&   1&   0&   0 \\
     0&  -1&   1&  -1&   1&  -1&   0&   1&   0&   0&  -1&  -2&   0&   0&   0&  -1&
  0&   0&   0&   0&   0&   0 \\
    -1&   1&   0&   0&   0&   0&  -1&   1&   1&   1&  -1&   0&   0&   0&   1&   1&
 -1&   0&  -1&   0&   0&   0 \\
     0&   0&   0&   0&   0&   0&   0&   1&   1&   1&  -1&   0&   0&   0&   1&   0&
 -1&   0&  -1&   0&   0&   0 \\
     0&   0&   0&   0&   1&  -1&   0&   0&   1&  -1&   0&  -1&  -1&   0&   0&  -1&
  0&  -1&   0&  -1&   1&   0 \\
     0&   0&  -1&   0&   0&   0&   1&   0&   0&   0&   0&   0&   0&   0&   1&  -1&
  0&   0&   1&   0&   0&   0 \\
\end{array}\right)$,

$s_{13}$=$\left(\begin{array}{cccccccccccccccccccccc}   1&   0&   0&
1& -1&   0&   1&  -1&   0&   0&   0&   1&   0& 0&   0&   0&
  0&   0&   0&  -1&   0&   0 \\
     0&   0&   0&   0&   0&   0&   1&  -1&  -1&  -1&   1&   0&   0&   1&   0&   0&
  1&  -1&   1&   0&   0&   0 \\
     0&   0&   0&   0&   1&  -1&   0&   0&   1&  -1&   0&  -1&  -1&   0&   0&  -1&
  0&  -1&   0&  -1&   1&   0 \\
    -1&   0&   0&   0&   1&   0&  -1&   1&   0&   0&   0&   0&   0&   0&   0&   0&
 -1&   1&   0&   1&   0&   0 \\
     0&   0&   0&   0&   1&  -1&   0&   0&   0&  -1&   0&   0&  -1&   0&  -1&  -1&
  0&   0&   0&  -1&   1&   0 \\
    -1&   0&   0&  -1&   1&  -1&   0&   0&   0&   0&   0&  -1&  -1&   0&   0&   0&
  1&  -1&   0&   0&   0&   0 \\
    -1&   0&   0&   0&   1&  -1&   0&   0&   0&   0&   1&   0&  -1&   0&  -1&  -1&
  1&   0&   1&   0&   1&   1 \\
    -1&   1&   0&   0&   0&   0&  -1&   0&   1&   1&   0&   1&  -1&   0&   0&   1&
  0&   0&  -1&   0&   0&   0 \\
\end{array}\right)$,

$s_{14}$=$\left(\begin{array}{cccccccccccccccccccccc} -1&   0& 0& 0&
0&  -1&   0&   1&   0&   1&   0&   0&   0&   0&   0&   0&
0&   0&   0&   0&   0&   1 \\
     0&  -1&   1&  -1&   1&  -1&   0&   1&   0&   0&  -1&  -2&   0&   0&   0&  -1&
  0&   0&   0&   0&   0&   0 \\
     0&   0&   0&   0&   0&   0&   0&   0&   1&   0&   0&   0&   0&   0&   0&   0&
  0&   0&   0&   0&   0&   0 \\
     0&   1&   0&   0&   0&   0&   0&   0&   1&   0&   0&   0&   0&   1&   1&   1&
  0&  -1&  -1&  -1&   0&  -1 \\
     0&   0&   0&   0&   0&   0&   0&   0&   1&   1&  -1&   0&   0&   0&   0&   1&
  0&   0&  -1&   0&   0&   0 \\
    -1&   0&   0&  -1&   1&  -1&   0&   0&   1&   1&   0&   0&  -1&   0&   0&   0&
  1&  -1&   0&   0&   1&   1 \\
    -1&   1&  -1&   0&   1&   0&   0&   0&   1&   0&   0&   0&   0&   1&   1&   0&
 -1&   0&  -1&   0&   1&  -1 \\
    -1&   1&  -1&   0&   0&  -1&   1&   0&   1&   1&   0&   1&  -1&   0&   1&   0&
  0&  -1&   0&  -1&   1&   0 \\
\end{array}\right)$,

$s_{15}$=$\left(\begin{array}{cccccccccccccccccccccc}  0&   0& 1& 0&
0&  -1&   0&   1&   1&   1&  -1&   0&  -1& -1&   0&   0&
  1&   0&   0&  -1&   0&   1 \\
    -1&   0&   0&   0&   0&  -1&   0&   1&   0&   1&   0&   0&   0&   0&   0&   0&
  0&   0&   0&   0&   0&   1 \\
     0&  -1&   0&  -1&   0&  -1&   1&   0&  -1&   0&   0&  -1&   0&   0&   0&  -1&
  1&  -1&   1&   0&   0&   0 \\
     0&   0&   0&  -1&   1&   0&   0&   0&   0&  -1&   0&  -1&   0&   1&   0&   0&
  0&  -1&  -1&   0&   0&  -1 \\
     0&  -1&   1&  -1&   1&  -1&  -1&   1&   0&   0&   0&  -1&  -1&  -1&  -1&  -1&
  1&   0&   0&   0&   0&   1 \\
    -1&   1&  -1&   0&   0&   0&   0&   0&   0&   0&   1&   0&   0&   1&   1&   0&
  0&  -1&   0&   0&   0&   0 \\
     0&  -1&   0&  -1&   0&   0&   0&   0&   0&   0&   0&  -1&   0&   0&   0&   0&
  0&  -1&   0&   0&   0&   0 \\
     0&   0&   0&   0&   0&   1&   0&   0&   0&   0&   0&   0&   0&   0&   0&   0&
  0&   0&   0&   0&   0&   0 \\
\end{array}\right)$,

$s_{16}$=$\left(\begin{array}{cccccccccccccccccccccc}  0&   1& 0& 1&
0&   0&   0&   1&   0&   0&   0&   0&   1& 1&   1&   0&
 -1&   1&   0&   0&   0&   0 \\
     0&   0&   0&   0&   0&   0&   0&   0&   0&   0&   1&   0&   0&   0&   0&   0&
  0&   0&   0&   0&   0&   0 \\
     0&   0&  -1&   0&   0&  -1&   1&  -1&   1&   0&   0&   0&  -1&   0&   0&  -1&
  0&  -1&   0&  -1&   1&   0 \\
     0&   1&  -1&   1&  -1&   0&   1&  -1&   1&   0&   0&   1&   0&   0&   1&   0&
  0&  -1&   0&  -1&   0&   0 \\
     0&   0&   0&   0&   0&  -1&   1&   0&   0&   0&   1&   0&   0&   0&   0&  -1&
  1&  -1&   1&  -1&   0&   1 \\
     0&   0&  -1&   0&   0&   1&   0&  -1&   0&   0&   1&   1&   0&   0&   0&   0&
  0&   0&   0&   0&   0&   0 \\
    -1&   0&  -1&   0&   0&   0&   1&  -1&   0&   0&   1&   1&   0&   1&   0&   0&
  0&  -1&   0&   0&   1&   0 \\
    -1&   1&  -1&   0&   0&   1&   0&   0&   0&   0&   0&   1&   0&   1&   1&   1&
 -1&   0&  -1&   0&   0&  -1 \\
\end{array}\right)$,

$s_{17}$=$\left(\begin{array}{cccccccccccccccccccccc}  1&   0& 0& 1&
0&   0&   0&   0&   0&  -1&   0&   0&   1&   0&   0&  -1&
-1&   1&   0&   0&   0&   0 \\
     0&   1&  -1&   1&  -1&   1&   0&   0&   0&   1&   0&   1&   1&   0&   1&   1&
 -1&   1&   0&   1&  -1&   0 \\
    -1&   1&   0&   0&   0&   0&   0&   0&   0&   1&   1&   1&   0&   0&   0&   1&
  1&   0&   0&   0&   0&   1 \\
     0&   0&   0&   0&   0&   0&   1&  -1&  -1&  -1&   1&   0&   0&   1&   0&   0&
  1&  -1&   1&   0&   0&   0 \\
     0&   0&   0&   1&   0&   0&   0&   0&  -1&   0&   1&   1&   1&   0&  -1&   0&
  0&   1&   1&   1&   0&   1 \\
     0&  -1&   1&  -1&   0&   0&   0&   1&  -1&   1&   0&  -1&   0&  -1&   0&   0&
  1&   0&   1&   1&  -1&   1 \\
     0&   1&   0&   1&  -1&   0&   0&   0&   1&   1&   0&   1&   0&   0&   0&   1&
  0&   0&   0&   0&   0&   1 \\
     1&   0&   1&   0&  -1&   0&   0&   1&   0&   0&  -1&  -1&   1&   0&   1&   0&
  0&   0&   0&   0&  -1&   0 \\
\end{array}\right)$,

$s_{18}$=$\left(\begin{array}{cccccccccccccccccccccc}  0&   0& 0& 0&
0&   0&   0&   0&   1&   0&   0&   0&   0& 0&   0&   0&
  0&   0&   0&   0&   0&   0 \\
    -1&   1&  -1&   0&   0&  -1&   1&   0&   1&   0&   0&   0&   0&   1&   1&   0&
  0&  -1&   0&  -1&   1&   0 \\
    -1&   0&   0&   0&   1&  -1&   0&   1&   0&   0&   0&   0&   0&   0&   0&  -1&
  0&   0&   0&   0&   1&   1 \\
     0&   0&   0&   0&   0&   0&   0&   0&   0&   0&   0&   0&   0&   0&   1&   0&
  0&   0&   0&   0&   0&   0 \\
     0&   0&   0&   0&   0&   0&   0&   0&   0&   0&   0&   0&   0&   0&   0&   0&
  1&   0&   0&   0&   0&   0 \\
     0&   0&   0&   0&   0&   0&   0&   0&  -1&  -1&   1&   0&   1&   1&   0&   0&
  0&   0&   0&   0&   0&   0 \\
     0&  -1&   0&  -1&   1&  -1&   1&   0&   0&   0&   0&  -1&  -1&   0&   0&  -1&
  1&  -1&   0&  -1&   1&   0 \\
     1&  -1&   1&  -1&   0&   1&   0&   0&  -1&  -1&   0&  -1&   0&   0&   0&   0&
  0&   0&   0&   0&  -1&  -1 \\
\end{array}\right)$,

$s_{19}$=$\left(\begin{array}{cccccccccccccccccccccc}  1&  -1& 1& 0&
0&   0&   0&   0&  -1&  -1&   0&  -1&   0& 0&  -1&  -1&
  0&   0&   0&   0&   0&   0 \\
     0&   0&  -1&   0&  -1&   1&   1&  -1&  -1&   0&   1&   1&   0&   0&   0&   0&
  0&   0&   1&   0&   0&   0 \\
     0&   0&   0&  -1&   0&   0&   0&   0&   0&   0&   0&   0&  -1&  -1&   0&   0&
  1&  -1&   0&  -1&   0&   0 \\
     0&   1&   0&   1&   0&   0&  -1&   1&   1&   1&   0&   1&   0&  -1&   0&   0&
  0&   1&   0&   0&   0&   1 \\
     1&   0&   0&   1&  -1&   1&   0&  -1&  -1&  -1&   1&   0&   1&   0&   0&   0&
  0&   0&   1&   0&  -1&   0 \\
     0&   0&   0&   0&  -1&   1&   0&   0&   0&   1&   0&   1&   0&  -1&   0&   1&
  0&   0&   0&   0&  -1&   0 \\
     0&   0&  -1&   0&   0&   1&   0&  -1&   0&   0&   1&   1&   0&   0&   0&   0&
  0&   0&   0&   0&   0&   0 \\
     0&   0&   1&   0&   1&   0&  -1&   1&   0&   0&   0&   0&   0&   0&  -1&   0&
  0&   1&  -1&   0&   0&   0 \\
\end{array}\right)$,

$s_{20}$=$\left(\begin{array}{cccccccccccccccccccccc}  0&   0& 0& 1&
-1&   1&   0&   0&  -1&   0&   1&   1&   1& 0&   0&   0&
  0&   1&   1&   1&  -1&   1 \\
     1&   0&   0&   1&  -1&   1&   0&   0&   0&   0&  -1&   0&   1&   0&   1&   0&
 -1&   1&   0&   0&  -1&  -1 \\
     0&   0&  -1&   0&   0&   0&   1&  -1&   0&  -1&   0&   0&   0&   1&   1&  -1&
 -1&  -1&   0&  -1&   1&  -1 \\
     0&   1&   0&   0&   0&   0&   0&   0&   0&   0&   0&   0&   0&   0&   0&   0&
  0&   0&   0&   0&   0&   0 \\
     1&  -1&   1&   0&   0&   0&   0&   0&  -1&   0&   0&   0&   0&  -1&  -1&  -1&
  1&   1&   1&   0&   0&   1 \\
     1&   0&   1&   0&   0&   0&  -1&   0&   1&   0&  -1&   0&  -1&  -1&  -1&   0&
  0&   0&  -1&  -1&   0&   0 \\
     0&   0&   0&   0&   0&   0&   0&   0&   0&   1&   0&   0&   0&   0&   0&   0&
  0&   0&   0&   0&   0&   0 \\
     0&   1&   0&   1&  -1&   0&   0&   0&   1&   1&   0&   1&   0&   0&   0&   1&
  0&   0&   0&   0&   0&   1 \\
\end{array}\right)$,

$s_{21}$=$\left(\begin{array}{cccccccccccccccccccccc} -1&   1& 0& 0&
1&  -1&  -1&   1&   1&   0&   0&   0&   0&   1&   0&   0&

 -1&   0&  -1&   0&   1&   0 \\
     0&  -1&   1&  -1&   1&   0&  -1&   0&  -1&  -1&   0&  -1&   0&   0&  -1&   0&
  0&   0&   0&   1&   0&   0 \\
     0&  -1&   0&   0&   1&   0&   0&  -1&   0&  -1&   0&   0&  -1&   0&  -1&  -1&
  0&   0&   0&   0&   1&   0 \\
     1&   0&   0&   1&  -1&   1&   0&   0&   0&   0&  -1&   0&   1&   0&   1&   0&
 -1&   1&   0&   0&  -1&  -1 \\
     0&   0&   0&   0&   0&   0&   0&   0&   0&   0&   0&   1&   0&   0&   0&   0&
  0&   0&   0&   0&   0&   0 \\
     1&  -1&   1&  -1&   0&   1&   0&   0&  -1&  -1&   0&  -1&   0&   0&   0&   0&
  0&   0&   0&   0&  -1&  -1 \\
     1&   0&   0&   0&  -1&   1&   0&  -1&  -1&  -1&   0&   0&   0&   0&   0&   0&
  0&   0&   0&   0&  -1&  -1 \\
     0&   1&   0&   0&  -1&   1&   0&   0&   0&   0&   0&   0&   1&   1&   1&   1&
 -1&   0&  -1&   0&  -1&  -1 \\
\end{array}\right)$,

$s_{22}$=$\left(\begin{array}{cccccccccccccccccccccc}-1&   0& 0& -1&
1&   0&  -1&   0&   0&   0&   0&   0&  -1& 0&  -1&   0&
  0&   0&  -1&   0&   1&   0 \\
     0&   0&   0&   0&   0&   1&   0&   0&   0&   0&   0&   0&   0&   0&   0&   0&
  0&   0&   0&   0&   0&   0 \\
     0&   0&   1&   0&   0&  -1&   0&   1&   1&   1&  -1&   0&  -1&  -1&   0&   0&
  1&   0&   0&  -1&   0&   1 \\
     0&   0&   0&   1&   0&  -1&   0&   0&   0&   0&   0&   0&   0&   0&   0&  -1&
  0&   0&   1&   0&   0&   1 \\
     0&  -1&   0&   0&   0&   0&   0&   0&  -1&  -1&   0&  -1&   0&   0&   0&  -1&
  0&   0&   1&   0&   0&   0 \\
     0&   0&  -1&   0&   0&   0&   1&   0&   0&   0&   0&   0&   0&   0&   1&  -1&
  0&   0&   1&   0&   0&   0 \\
     1&  -1&   0&   1&   0&   0&   0&   0&   0&   0&   0&   0&   0&  -1&  -1&  -1&
  0&   1&   1&   0&   0&   1 \\
     0&   0&   0&   0&   1&   0&  -1&   0&   0&   0&   0&   0&   0&   0&  -1&   0&
  0&   1&   0&   1&   0&   0 \\
\end{array}\right)$,

$s_{23}$=$\left(\begin{array}{cccccccccccccccccccccc}  1&  -1& 1& 0&
0&   0&   0&   1&  -1&   0&  -1&  -1&   1& 0&   0&   0&
  0&   1&   0&   1&  -1&   0 \\
     0&   0&   0&   0&   0&   0&   0&   0&   0&   0&   0&   0&   0&   1&   0&   0&
  0&   0&   0&   0&   0&   0 \\
     0&   0&   0&   0&   0&   0&   0&  -1&   0&  -1&   1&   0&  -1&   0&  -1&   0&
  1&  -1&   0&  -1&   0&   0 \\
     0&   0&   0&   0&   1&  -1&   0&   0&   1&  -1&   0&  -1&  -1&   0&   0&  -1&
  0&  -1&   0&  -1&   1&   0 \\
     1&  -1&   1&   0&   0&   0&   0&   0&  -1&  -1&   0&  -1&   0&   0&  -1&  -1&
  0&   0&   0&   0&   0&   0 \\
     0&   0&   0&   0&   0&  -1&   1&   0&   0&   0&   1&   0&   0&   0&   0&  -1&
  1&  -1&   1&  -1&   0&   1 \\
     1&   0&   1&   0&   0&   0&  -1&   0&   1&   0&  -1&   0&  -1&  -1&  -1&   0&
  0&   0&  -1&  -1&   0&   0 \\
     0&   0&   0&   0&   0&  -1&   0&   1&   1&   1&  -1&   0&   0&  -1&   0&   0&
  0&   0&   0&   0&   0&   1 \\
\end{array}\right)$,

$s_{24}$=$\left(\begin{array}{cccccccccccccccccccccc}-1&   1& -1& 0&
0&   0&   0&   0&   0&   0&   1&   0&   0& 1&   1&   0&
  0&  -1&   0&   0&   0&   0 \\
     0&   0&  -1&   0&   0&   0&   1&   0&   0&   0&   0&   0&   0&   0&   1&  -1&
  0&   0&   1&   0&   0&   0 \\
     1&  -1&   1&   0&   0&   0&   0&   0&  -1&   0&   0&   0&   0&  -1&  -1&  -1&
  1&   1&   1&   0&   0&   1 \\
     1&  -1&   1&   0&   0&  -1&   0&   1&   0&   0&  -1&  -1&   0&  -1&   0&  -1&
  0&   0&   0&  -1&   0&   0 \\
     1&   0&   0&   0&  -1&   1&   1&   0&  -1&   0&   0&   0&   1&   0&   1&   0&
  0&   0&   0&   0&  -1&  -1 \\
     0&   0&   0&   0&  -1&   0&   1&   0&  -1&   0&   0&   0&   1&   0&   1&   0&
  0&   0&   1&   0&  -1&   0 \\
     0&  -1&   0&   0&   0&  -1&   1&   0&   0&   1&   0&   0&   0&  -1&   0&  -1&
  1&   0&   1&   0&   0&   1 \\
    -1&   0&   0&   0&   0&  -1&   0&   0&   0&   0&   1&   0&   0&   0&   0&   0&
  1&  -1&   1&   0&   0&   1 \\
\end{array}\right)$,

 Obviously,  $s_1$ is the unity element
and $G= G^{s_1}$.

The character table of $G^{s_1} =G$:\\
\begin{tabular}{c|cccccccccccccccccccc}
  & & & & & & & & & &10 & & & & & & & & & &20\\\hline\\$\chi_{1}^{(1)}$&1&1&1&1&1&1&1&1&1&1&1&1&1&1&1&1&1&1&1&1
\\$\chi_{1}^{(2)}$&22&-2&-2&.&.&-3&-5&1&6&3&-1&-1&1&1&-1&2&4&1&1&.
\\$\chi_{1}^{(3)}$&231&2&2&.&.&6&15&2&7&7&.&.&.&.&-1&-1&6&.&.&.
\\$\chi_{1}^{(4)}$&252&1&1&-1&-1&2&9&-2&28&1&.&.&.&.&1&4&9&.&.&-1
\\$\chi_{1}^{(5)}$&770&A&/A&A&/A&-5&-13&1&-14&7&.&.&.&.&1&-2&5&-1&-1&.
\\$\chi_{1}^{(6)}$&770&/A&A&/A&A&-5&-13&1&-14&7&.&.&.&.&1&-2&5&-1&-1&.
\\$\chi_{1}^{(7)}$&896&.&.&2&2&-4&32&.&.&.&.&.&.&.&.&.&-4&-1&-1&D
\\$\chi_{1}^{(8)}$&896&.&.&2&2&-4&32&.&.&.&.&.&.&.&.&.&-4&-1&-1&/D
\\$\chi_{1}^{(9)}$&1750&.&.&.&.&.&-5&.&70&-5&.&.&.&.&-1&2&13&-2&-2&1
\\$\chi_{1}^{(10)}$&3520&-1&-1&1&1&-5&-44&-1&64&4&1&1&-1&-1&.&.&10&1&1&.
\\$\chi_{1}^{(11)}$&3520&1&1&1&1&-5&-44&1&-64&-4&-1&-1&-1&-1&.&.&10&1&1&.
\\$\chi_{1}^{(12)}$&4500&.&.&.&.&.&45&.&20&5&-1&-1&-1&-1&1&4&-9&.&.&1
\\$\chi_{1}^{(13)}$&4752&-1&-1&-1&-1&2&54&2&-48&-6&1&1&-1&-1&.&.&.&.&.&.
\\$\chi_{1}^{(14)}$&5103&.&.&.&.&3&.&3&63&.&.&.&.&.&.&3&.&.&.&-1
\\$\chi_{1}^{(15)}$&5544&-1&-1&1&1&19&36&-1&-56&4&.&.&.&.&.&.&9&.&.&.
\\$\chi_{1}^{(16)}$&8019&.&.&.&.&-6&.&.&-45&.&B&/B&B&/B&.&3&.&.&.&.
\\$\chi_{1}^{(17)}$&8019&.&.&.&.&-6&.&.&-45&.&/B&B&/B&B&.&3&.&.&.&.
\\$\chi_{1}^{(18)}$&8250&.&.&.&.&.&15&.&10&-5&-/B&-B&/B&B&1&-2&6&.&.&.
\\$\chi_{1}^{(19)}$&8250&.&.&.&.&.&15&.&10&-5&-B&-/B&B&/B&1&-2&6&.&.&.
\\$\chi_{1}^{(20)}$&9625&.&.&.&.&.&40&.&105&.&.&.&.&.&.&-3&-5&1&1&.
\\$\chi_{1}^{(21)}$&9856&.&.&.&.&6&-80&.&.&.&.&.&.&.&.&.&-8&C&/C&.
\\$\chi_{1}^{(22)}$&9856&.&.&.&.&6&-80&.&.&.&.&.&.&.&.&.&-8&/C&C&.
\\$\chi_{1}^{(23)}$&10395&-/A&-A&/A&A&-5&27&-1&-21&3&.&.&.&.&-1&-1&.&.&.&.
\\$\chi_{1}^{(24)}$&10395&-A&-/A&A&/A&-5&27&-1&-21&3&.&.&.&.&-1&-1&.&.&.&.
\end{tabular}

\begin{tabular}{c|cccc}
  & & & &\\\hline
\\$\chi_{1}^{(1)}$&1&1&1&1
\\$\chi_{1}^{(2)}$&.&.&.&2
\\$\chi_{1}^{(3)}$&.&-1&-2&1
\\$\chi_{1}^{(4)}$&-1&.&1&2
\\$\chi_{1}^{(5)}$&.&.&1&.
\\$\chi_{1}^{(6)}$&.&.&1&.
\\$\chi_{1}^{(7)}$&/D&.&.&1
\\$\chi_{1}^{(8)}$&D&.&.&1
\\$\chi_{1}^{(9)}$&1&.&1&.
\\$\chi_{1}^{(10)}$&.&.&-2&.
\\$\chi_{1}^{(11)}$&.&.&2&.
\\$\chi_{1}^{(12)}$&1&.&-1&.
\\$\chi_{1}^{(13)}$&.&.&.&2
\\$\chi_{1}^{(14)}$&-1&1&.&-2
\\$\chi_{1}^{(15)}$&.&.&1&-1
\\$\chi_{1}^{(16)}$&.&-1&.&-1
\\$\chi_{1}^{(17)}$&.&-1&.&-1
\\$\chi_{1}^{(18)}$&.&.&-2&.
\\$\chi_{1}^{(19)}$&.&.&-2&.
\\$\chi_{1}^{(20)}$&.&-1&3&.
\\$\chi_{1}^{(21)}$&.&.&.&1
\\$\chi_{1}^{(22)}$&.&.&.&1
\\$\chi_{1}^{(23)}$&.&1&.&.
\\$\chi_{1}^{(24)}$&.&1&.&.

\end{tabular}

 \noindent  \noindent where
A = -E(15)$^7$-E(15)$^{11}$-E(15)$^{13}$-E(15)$^{14}$
  = (-1+ER(-15))/2 = b15;
B = -E(7)-E(7)$^2$-E(7)$^4$
  = (1-ER(-7))/2 = -b7;
C = 2*E(3)-E(3)$^2$
  = (-1+3*ER(-3))/2 = 1+3b3;
D = E(11)+E(11)$^3$+E(11)$^4$+E(11)$^5$+E(11)$^9$
  = (-1+ER(-11))/2 = b11;

The generators of $G^{s_2}$ are:\\
$\left(\begin{array}{cccccccccccccccccccccc}0& 0& 0&
0& 0& 0& 0& -1& 0& -1& 0& -1& 0& 1& 0& 0& 0& -1& 0& 0& 0& -1 \\
       0& 0& 0& 1& 0& -1& 0& 0& 0& 0& 0& 0& 0& 0& 0& -1& 0& 0& 1& 0& 0& 1 \\
       0& 1& 1& 0& 0& 0& -1& 1& 0& 0& 0& 0& 0& 0& 0& 1& 0& 0& -1& 0& -1& 0 \\
       0& 0& 1& -1& 1& -1& 0& 1& 0& 0& 0& -1& 0& 0& 0& 0& 1& 0& 0& 0& 0& 0 \\
       0& 0& 0& 0& 0& 0& 1& -1& -1& -1& 1& 0& 0& 1& 0& 0& 1& -1& 1& 0& 0& 0 \\
       0& 0& 0& 1& 0& 0& 0& 0& -1& 0& 1& 1& 1& 0& -1& 0& 0& 1& 1& 1& 0& 1 \\
       0& 0& 0& 0& -1& 0& 1& 0& -1& 0& 0& 0& 1& 0& 1& 0& 0& 0& 1& 0& -1& 0 \\
       -1& 0& 0& -1& 0& 0& 0& 0& 0& 1& 0& 0& 0& 0& 0& 1& 0& 0& 0& 1& 0& 0 \\
       0& 0& 0& 0& -1& 1& 0& 0& 0& 1& 0& 1& 0& -1& 0& 1& 0& 0& 0& 0& -1& 0 \\
       1& 0& 1& 1& 0& 0& -1& 0& 0& -1& 0& 0& 0& 0& -1& 0& 0& 1& 0& 0& 0& 0 \\
       0& 0& 1& -1& 0& 0& -1& 1& 0& 0& -1& -1& 0& 0& 0& 1& 0& 0& -1& 0& -1& -1 \\
       0& 0& 0& -1& 0& 1& 0& -1& 0& 0& 0& 0& -1& 0& 0& 1& 0& -1& -1& 0& 0& -1 \\
       -1& 0& 0& 0& 1& 0& -1& 0& 0& 0& 1& 1& -1& 0& -1& 0& 0& 0& 0& 0& 1& 1 \\
       1& 0& 0& 1& -1& 1& 0& -1& 0& 0& 0& 1& 0& -1& -1& 0& 0& 1& 0& 0& 0& 0 \\
       0& 0& 0& -1& 1& 0& 0& 0& 0& -1& 0& -1& 0& 1& 0& 0& 0& -1& -1& 0& 0& -1 \\
       -1& 1& -1& 0& 0& -1& 1& 0& 1& 0& 0& 0& 0& 1& 1& 0& 0& -1& 0& -1& 1& 0 \\
       0& -1& 0& -1& 1& 0& 0& 0& 0& 0& 0& 0& -1& -1& -1& -1& 0& 0& 0& 0& 1& 0 \\
       1& -1& 0& 0& -1& 0& 1& 0& 0& 0& -1& 0& 0& -1& 0& -1& 0& 0& 0& -1& 0& 0 \\
       0& 1& 0& 1& -1& 0& 0& 0& 1& 1& 0& 1& 0& 0& 0& 1& 0& 0& 0& 0& 0& 1 \\
       0& 0& 0& -1& 0& 0& 0& 0& 0& 0& 0& -1& 0& 0& 1& 0& 0& -1& 0& 0& -1& -1 \\
       -1& 1& 0& 0& 1& -1& -1& 1& 1& 0& 0& 0& 0& 1& 0& 0& -1& 0& -1& 0& 1& 0 \\
       0& 0& -1& 1& 0& 0& 1& -1& -1& -1& 1& 0& 1& 1& 0& -1& 0& 0& 1& 0& 0& 0
       \\
\end{array}\right)$.

The representatives of conjugacy classes of   $G^{s_2}$ are:\\
$\left(\begin{array}{cccccccccccccccccccccc}1& 0& 0& 0& 0& 0& 0& 0& 0& 0& 0& 0& 0& 0& 0& 0& 0& 0& 0& 0& 0& 0 \\
       0& 1& 0& 0& 0& 0& 0& 0& 0& 0& 0& 0& 0& 0& 0& 0& 0& 0& 0& 0& 0& 0 \\
       0& 0& 1& 0& 0& 0& 0& 0& 0& 0& 0& 0& 0& 0& 0& 0& 0& 0& 0& 0& 0& 0 \\
       0& 0& 0& 1& 0& 0& 0& 0& 0& 0& 0& 0& 0& 0& 0& 0& 0& 0& 0& 0& 0& 0 \\
       0& 0& 0& 0& 1& 0& 0& 0& 0& 0& 0& 0& 0& 0& 0& 0& 0& 0& 0& 0& 0& 0 \\
       0& 0& 0& 0& 0& 1& 0& 0& 0& 0& 0& 0& 0& 0& 0& 0& 0& 0& 0& 0& 0& 0 \\
       0& 0& 0& 0& 0& 0& 1& 0& 0& 0& 0& 0& 0& 0& 0& 0& 0& 0& 0& 0& 0& 0 \\
       0& 0& 0& 0& 0& 0& 0& 1& 0& 0& 0& 0& 0& 0& 0& 0& 0& 0& 0& 0& 0& 0 \\
       0& 0& 0& 0& 0& 0& 0& 0& 1& 0& 0& 0& 0& 0& 0& 0& 0& 0& 0& 0& 0& 0 \\
       0& 0& 0& 0& 0& 0& 0& 0& 0& 1& 0& 0& 0& 0& 0& 0& 0& 0& 0& 0& 0& 0 \\
       0& 0& 0& 0& 0& 0& 0& 0& 0& 0& 1& 0& 0& 0& 0& 0& 0& 0& 0& 0& 0& 0 \\
       0& 0& 0& 0& 0& 0& 0& 0& 0& 0& 0& 1& 0& 0& 0& 0& 0& 0& 0& 0& 0& 0 \\
       0& 0& 0& 0& 0& 0& 0& 0& 0& 0& 0& 0& 1& 0& 0& 0& 0& 0& 0& 0& 0& 0 \\
       0& 0& 0& 0& 0& 0& 0& 0& 0& 0& 0& 0& 0& 1& 0& 0& 0& 0& 0& 0& 0& 0 \\
       0& 0& 0& 0& 0& 0& 0& 0& 0& 0& 0& 0& 0& 0& 1& 0& 0& 0& 0& 0& 0& 0 \\
       0& 0& 0& 0& 0& 0& 0& 0& 0& 0& 0& 0& 0& 0& 0& 1& 0& 0& 0& 0& 0& 0 \\
       0& 0& 0& 0& 0& 0& 0& 0& 0& 0& 0& 0& 0& 0& 0& 0& 1& 0& 0& 0& 0& 0 \\
       0& 0& 0& 0& 0& 0& 0& 0& 0& 0& 0& 0& 0& 0& 0& 0& 0& 1& 0& 0& 0& 0 \\
       0& 0& 0& 0& 0& 0& 0& 0& 0& 0& 0& 0& 0& 0& 0& 0& 0& 0& 1& 0& 0& 0 \\
       0& 0& 0& 0& 0& 0& 0& 0& 0& 0& 0& 0& 0& 0& 0& 0& 0& 0& 0& 1& 0& 0 \\
       0& 0& 0& 0& 0& 0& 0& 0& 0& 0& 0& 0& 0& 0& 0& 0& 0& 0& 0& 0& 1& 0 \\
       0& 0& 0& 0& 0& 0& 0& 0& 0& 0& 0& 0& 0& 0& 0& 0& 0& 0& 0& 0& 0& 1 \\
\end{array}\right)$,

$\left(\begin{array}{cccccccccccccccccccccc}
 -1& 0& -1& 0& 0& 0& 1& -1& 0& 0& 1& 0& 0& 1& 1& 0& 0& -1& 1& 0& 0& 0 \\
       0& 0& 0& 0& -1& 0& 0& 0& 0& 1& 0& 0& 0& -1& 0& 0& 1& 0& 1& 0& -1& 1 \\
       0& 0& 1& 0& 1& -1& -1& 1& 0& 0& 0& -1& 0& 0& -1& 0& 1& 0& 0& 0& 0& 1 \\
       0& 0& 0& 1& 0& 0& 0& 0& -1& 0& 1& 1& 1& 0& -1& 0& 0& 1& 1& 1& 0& 1 \\
       0& -1& 1& -1& 0& 0& 0& 1& -1& 1& 0& -1& 0& -1& 0& 0& 1& 0& 1& 1& -1& 1 \\
       0& 0& 0& 0& -1& 1& 0& 0& -1& 0& 0& 0& 1& 0& 0& 1& 0& 0& 0& 1& -1& 0 \\
       0& 0& 1& 0& 0& 0& -1& 1& -1& 0& 0& -1& 1& 0& 0& 0& 0& 1& 0& 1& -1& 0 \\
       -1& 1& -1& 1& 0& 0& 0& -1& 0& -1& 1& 1& 0& 1& 0& 0& -1& 0& 0& 0& 1& 0 \\
       1& 0& 1& 0& 0& 0& 0& 0& 0& -1& -1& -1& 0& 0& 0& 0& 0& 0& -1& -1& 0& -1 \\
       0& 0& -1& 1& -1& 0& 1& 0& 0& 1& 0& 1& 1& 0& 1& 0& 0& 0& 1& 0& 0& 1 \\
       -1& 0& 0& 0& 1& -1& 0& 0& 0& 0& 0& 0& 0& 1& 0& 0& 0& 0& 0& 0& 1& 0 \\
       0& 0& 0& 0& 1& -1& 0& 0& 1& -1& 0& -1& -1& 0& 0& -1& 0& -1& 0& -1& 1& 0 \\
       0& 0& 0& -1& 0& 0& 0& 0& 0& 0& 0& 0& -1& -1& 0& 0& 1& -1& 0& -1& 0& 0 \\
       1& 0& -1& 1& -1& 1& 1& -1& 0& -1& 0& 0& 1& 1& 1& 0& -1& 0& 0& 0& 0& -1 \\
       0& -1& 1& -1& 1& -1& 0& 0& 0& 0& 0& -1& -1& 0& -1& 0& 1& -1& 0& 0& 0& 0 \\
       -1& 0& 0& -1& 1& 0& 0& 0& -1& 0& 1& 0& 0& 1& 0& 0& 0& 0& 0& 1& 0& 0 \\
       1& 0& 0& 0& -1& 0& 1& -1& 0& -1& 0& 0& 0& 0& 0& 0& 0& -1& 0& -1& 0& -1 \\
       1& -1& 1& -1& 0& 0& 0& 0& 0& 0& -1& -1& 0& 0& 0& 0& 0& 0& -1& 0& 0& -1 \\
       0& 0& 0& 0& 0& 0& 0& 0& 0& 0& 0& 0& 0& 0& 0& 0& 0& 0& 0& 0& 0& 1 \\
       0& 1& 0& 1& 0& -1& 0& 0& 1& 0& 0& 0& 0& 0& 0& 0& 0& 0& 0& -1& 0& 0 \\
       -1& 1& -1& 0& 0& 0& 0& -1& 0& 0& 1& 1& 0& 1& 0& 1& 0& -1& 0& 0& 0& 0 \\
       0& -1& 0& -1& 0& 1& 0& 0& -1& 0& 0& -1& 0& 0& 0& 0& 0& 0& 0& 1& -1& -1
       \\
\end{array}\right)$,

$\left(\begin{array}{cccccccccccccccccccccc}
-1& 0& -1& 0& 0& 1& 0& -1& -1& 0& 1& 1& 0& 0& 0& 0& 0& 0& 1& 1& 0& 0 \\
       0& 0& 0& 0& 0& 0& 1& -1& -1& -1& 1& 0& 0& 1& 0& 0& 1& -1& 1& 0& 0& 0 \\
       1& 0& 0& 1& 0& -1& 1& 0& 1& 0& -1& 0& 0& 0& 0& -1& 0& 0& 0& -1& 1& 0 \\
       1& 0& 1& 0& 0& 0& -1& 0& 1& 0& -1& 0& -1& -1& -1& 0& 0& 0& -1& -1& 0& 0 \\
       0& 1& -1& 1& -1& 0& 1& -1& 1& 0& 0& 1& 0& 0& 1& 0& 0& -1& 0& -1& 0& 0 \\
       1& -1& 0& 0& 0& -1& 1& -1& 0& -1& 0& -1& 0& 0& 0& -1& 1& -1& 1& -1& 0& 0 \\
       0& 0& 1& 0& 0& -1& 0& 1& 1& 1& -1& 0& -1& -1& 0& 0& 1& 0& 0& -1& 0& 1 \\
       0& -1& 1& -1& 1& 0& -1& 0& 0& 0& 0& -1& -1& -1& -1& 0& 1& 0& 0& 0& 0& 0 \\
       0& 0& 0& 0& 0& -1& 0& 1& 1& 1& -1& -1& 0& 0& 1& 0& 0& 0& 0& 0& 0& 0 \\
       0& 0& -1& 1& 0& 0& 0& -1& 0& -1& 1& 1& 0& 0& -1& -1& 0& 0& 1& 0& 1& 1 \\
       0& -1& 1& 0& 0& 0& -1& 1& 0& 1& -1& 0& 0& -1& -1& 0& 0& 1& 0& 1& 0& 1 \\
       -1& 1& -1& 0& 1& 0& 0& 0& 1& 0& 0& 0& 0& 1& 1& 0& -1& 0& -1& 0& 1& -1 \\
       1& 0& -1& 1& -1& 1& 1& -1& 0& -1& 0& 0& 1& 1& 1& 0& -1& 0& 0& 0& 0& -1 \\
       -1& 0& -1& 0& 1& -1& 0& 0& 0& 0& 1& 0& 0& 0& 0& -1& 0& 0& 1& 0& 1& 1 \\
       0& 0& 0& 0& 0& 0& 0& 0& 0& 0& -1& 0& 0& 0& 0& 0& -1& 0& -1& 0& 0& -1 \\
       0& 0& 0& 0& -1& 1& 0& 0& 0& 1& 0& 1& 0& -1& 0& 1& 0& 0& 0& 0& -1& 0 \\
       0& 0& 0& 0& 1& 0& -1& 0& 0& -1& 0& -1& 0& 1& 0& 0& -1& 0& -1& 0& 0& -1 \\
       -1& 0& 0& 0& 0& -1& 0& 1& 0& 1& 0& 0& 0& 0& 0& 0& 0& 0& 0& 0& 0& 1 \\
       0& 0& -1& 0& 0& 0& 1& -1& 0& 0& 1& 0& 0& 0& 0& 0& 1& -1& 1& 0& 0& 0 \\
       0& 0& 1& -1& 1& 0& -1& 1& 0& 0& -1& -1& -1& 0& 0& 0& 0& 0& -1& 0& 0& -1 \\
       1& -1& 0& 0& 0& 1& 0& -1& -1& -1& 0& 0& 0& 0& -1& 0& 0& 0& 0& 0& 0& -1 \\
       0& 0& 0& 0& -1& 0& 1& 0& -1& 0& 0& 0& 1& 0& 1& 0& 0& 0& 1& 0& -1& 0
       \\
\end{array}\right)$,

$\left(\begin{array}{cccccccccccccccccccccc}
-1& 0& 0& 0& 1& -1& 0& 0& 0& 0& 0& 0& 0& 1& 0& 0& 0& 0& 0& 0& 1& 0 \\
       -1& 1& 0& 0& 1& 0& -1& 0& 1& 0& 0& 1& -1& 0& -1& 1& 0& 0& -1& 0& 1& 0 \\
       1& 0& 0& 0& 0& 0& 0& 0& 1& 0& -1& 0& -1& -1& 0& 0& 0& 0& -1& -1& 0& -1 \\
       0& 0& 0& 0& 0& 0& 0& 1& 1& 1& -1& 0& 0& 0& 1& 0& -1& 0& -1& 0& 0& 0 \\
       0& 0& -1& 1& 0& 0& 0& -1& 1& 0& 0& 1& 0& 0& 0& 0& -1& 0& 0& 0& 1& 0 \\
       1& 0& 1& 0& 0& 0& -1& 0& 1& 0& -1& 0& -1& -1& -1& 0& 0& 0& -1& -1& 0& 0 \\
       0& 0& 0& 0& 0& 1& -1& 0& 0& 0& 0& 0& 0& 0& 0& 1& 0& 0& 0& 1& -1& 0 \\
       0& -1& 1& -1& 0& 0& 0& 0& -1& 0& 0& -1& 0& 0& 0& 0& 1& -1& 0& 0& -1& 0 \\
       1& 0& 0& 0& -1& 1& 0& -1& -1& -1& 0& 0& 0& 0& 0& 0& 0& 0& 0& 0& -1& -1 \\
       0& 0& 1& 0& 1& -1& -1& 1& 1& 1& -1& 0& -1& -1& -1& 0& 0& 1& -1& 0& 1& 1 \\
       1& -1& 1& 0& 0& 0& 0& 1& -1& 0& -1& -1& 1& 0& 0& 0& 0& 1& 0& 1& -1& 0 \\
       0& 0& -1& 0& -1& 0& 1& -1& 0& 0& 0& 0& 0& 0& 1& 0& 0& -1& 0& -1& 0& -1 \\
       1& 0& 0& 1& -1& 1& 0& -1& 0& 0& 0& 1& 0& -1& -1& 0& 0& 1& 0& 0& 0& 0 \\
       0& 0& 1& -1& 0& -1& 0& 0& 0& 0& 0& -1& -1& 0& 0& 0& 1& -1& 0& -1& 0& 0 \\
       0& 0& -1& 0& 0& 0& 1& 0& 0& 0& 0& 0& 1& 1& 1& 0& -1& 0& 0& 0& 0& -1 \\
       0& -1& 0& 0& 1& 0& 0& 0& -1& -1& 1& -1& 1& 1& 0& -1& 0& 0& 1& 1& 0& 0 \\
       0& 1& 0& 0& -1& 1& 0& 0& 0& 0& 0& 0& 1& 1& 1& 1& -1& 0& -1& 0& -1& -1 \\
       0& 0& 0& 0& 0& 0& 0& 0& -1& -1& 1& -1& 1& 1& 0& 0& 0& 0& 1& 1& -1& 0 \\
       0& 0& 0& 0& 1& -1& 0& -1& 0& -1& 1& 0& -1& 0& -1& -1& 1& -1& 0& -1& 1& 0 \\
       -1& 1& -1& 0& 0& 1& 0& 0& 0& 0& 0& 1& 0& 1& 1& 1& -1& 0& -1& 0& 0& -1 \\
       0& 0& 1& 0& 1& 0& -1& 1& 0& 0& 0& 0& 0& 0& -1& 0& 0& 1& -1& 0& 0& 0 \\
       0& 0& 0& 0& 0& 0& 0& 0& 0& 0& 0& 0& 0& 0& 0& 0& 0& 0& 1& 0& 0& 0
       \\
\end{array}\right)$,

$\left(\begin{array}{cccccccccccccccccccccc}
0& -1& 0& -1& 0& 0& 1& -1& -1& -1& 1& -1& 0& 0& 0& -1& 1& -1& 1& 0& 0& 0 \\
       0& -1& 0& 0& 0& 1& 0& 0& -1& 0& 0& 0& 1& 0& 0& 0& -1& 1& 0& 1& 0& 0 \\
       1& 0& 0& 1& -1& 0& 0& 0& 0& 0& -1& 0& 1& 0& 0& 0& -1& 1& 0& 0& 0& 0 \\
       0& 0& 0& 0& 0& 0& 0& 0& 0& 0& 0& 0& 0& 1& 0& 0& 0& 0& 0& 0& 0& 0 \\
       0& 0& 0& 0& 0& 0& 0& 0& 0& 0& 0& 0& 0& 0& 0& 0& 0& 1& 0& 0& 0& 0 \\
       0& -1& 0& -1& 0& 0& 0& 0& 0& 1& -1& 0& -1& -1& 0& 0& 0& 0& 0& 0& 0& 0 \\
       0& 0& 0& 0& 0& 0& 0& 0& 1& 0& -1& 0& -1& 0& 0& 0& -1& 0& -1& -1& 1& -1 \\
       1& -1& 0& 0& -1& 0& 1& -1& 0& 0& 0& 0& -1& -1& 0& -1& 1& -1& 1& -1& 0& 0 \\
       0& 0& 0& 0& 0& 0& -1& 0& 0& 0& 0& 0& -1& -1& -1& 0& 0& 0& 0& 0& 0& 0 \\
       0& 0& 0& -1& 0& 1& 0& 0& -1& 0& 0& -1& 1& 1& 1& 1& 0& 0& 0& 1& -1& -1 \\
       0& 0& 0& 0& 0& 0& 0& -1& 0& -1& 0& -1& 0& 1& 0& 0& 0& -1& 0& 0& 0& -1 \\
       1& 0& 0& 1& -1& -1& 1& 0& 0& 0& 0& 0& 0& -1& 0& -1& 1& 0& 1& -1& 0& 1 \\
       0& 1& -1& 1& -1& 1& 0& 0& 0& 1& 0& 1& 1& 0& 1& 1& -1& 1& 0& 1& -1& 0 \\
       0& 0& 0& -1& 0& 0& 0& 0& 0& 0& 0& 0& -1& -1& 0& 0& 1& -1& 0& -1& 0& 0 \\
       0& 0& 1& 0& 0& -1& 0& 0& 0& 0& 0& 0& 0& 0& -1& 0& 1& 0& 0& 0& 0& 1 \\
       0& 0& 0& 0& 0& 0& 0& 0& 1& 0& 0& 0& 0& 0& 0& 0& 0& 0& 0& 0& 0& 0 \\
       0& 0& 0& 0& 0& 0& 0& 0& 0& 1& 0& 1& -1& -1& -1& 0& 1& 0& 0& 0& 0& 1 \\
       -1& 1& 0& 0& 1& 0& -1& 0& 1& 0& 0& 1& -1& 0& -1& 1& 0& 0& -1& 0& 1& 0 \\
       1& -1& 0& 0& -1& 1& 0& 0& -1& 0& 0& 0& 1& -1& 0& 0& 0& 1& 1& 1& -1& 0 \\
       1& -1& 1& 0& 0& 0& 0& 0& -1& -1& 0& -1& 0& 0& -1& -1& 0& 0& 0& 0& 0& 0 \\
       1& -1& 0& 0& -1& 1& 1& -1& -1& 0& 0& 0& 1& 0& 0& 0& 0& 0& 1& 1& -1& 0 \\
       -1& 0& 0& -1& 1& -1& 0& 0& 1& 0& 0& 0& -1& 0& 0& 0& 0& -1& -1& -1& 1& 0
       \\

\end{array}\right)$,

$\left(\begin{array}{cccccccccccccccccccccc}
0& -1& 0& 0& 0& -1& 1& -1& 0& 0& 0& 0& -1& -1& -1& -1& 1& -1& 1& -1& 1& 1 \\
       1& -1& 1& 0& -1& 0& 1& 0& -1& 0& 0& -1& 1& 0& 0& 0& 1& 0& 1& 0& -1& 0 \\
       1& 0& 1& 0& -1& 0& 0& 1& 0& 1& -1& 0& 0& -1& 0& 0& 0& 1& 0& 0& -1& 0 \\
       0& 1& 0& 1& 0& 0& -1& 1& 1& 0& -1& 0& 0& 0& 0& 0& -1& 1& -1& 0& 0& 0 \\
       0& 0& 0& 0& -1& 0& 0& 0& 0& 1& 0& 0& 0& -1& 0& 0& 1& 0& 1& 0& -1& 1 \\
       0& 0& 0& 0& 0& 0& 0& 0& 0& 0& 0& 0& 0& 1& 0& 0& 0& 0& 0& 0& 0& 0 \\
       0& -1& 0& -1& 0& 1& 0& 0& -1& 0& 0& -1& 0& 0& 0& 0& 0& 0& 0& 1& -1& -1 \\
       0& -1& 0& 0& 1& 0& 0& -1& 0& -1& 0& 0& -1& 0& -1& -1& 0& 0& 0& 0& 1& 0 \\
       -1& 0& 0& -1& 1& 0& 0& 0& -1& 0& 1& 0& 0& 1& 0& 0& 0& 0& 0& 1& 0& 0 \\
       0& 1& 0& 1& -1& 0& 0& 0& 0& 0& 0& 1& 0& 0& 0& 0& 0& 0& 0& -1& 0& 0 \\
       0& 0& 0& 0& 0& 0& 0& 0& 0& 0& -1& 0& 0& 0& 0& 0& -1& 0& -1& 0& 0& -1 \\
       0& 0& 0& 0& 0& 0& 0& -1& 0& 0& 1& 1& -1& -1& -1& 0& 1& 0& 1& 0& 0& 1 \\
       0& 1& -1& 1& -1& 0& 1& -1& 0& 0& 1& 1& 1& 1& 1& 0& 0& 0& 1& 0& 0& 0 \\
       -1& 0& 0& 0& 1& 0& -1& 0& 0& 0& 1& 1& -1& 0& -1& 0& 0& 0& 0& 0& 1& 1 \\
       0& 0& 0& 0& 0& -1& 0& 1& 1& 1& -1& 0& 0& -1& 0& 0& 0& 0& 0& 0& 0& 1 \\
       1& -1& 0& 0& 0& 0& 0& 0& 0& 0& -1& -1& 0& -1& 0& -1& 0& 0& 0& 0& 0& 0 \\
       -1& 1& -1& 1& 1& -1& 0& 0& 1& 0& 1& 1& 0& 1& 0& 0& 0& 0& 0& 0& 1& 1 \\
       -1& 0& 0& -1& 1& 0& -1& 1& 0& 1& 0& 0& 0& 0& 0& 1& 0& 0& -1& 1& 0& 0 \\
       1& -1& 1& 0& 0& 0& 0& 0& -1& 0& 0& 0& 0& -1& -1& -1& 1& 1& 1& 0& 0& 1 \\
       0& 0& 0& 0& 0& 0& 1& 0& 0& 0& 0& 0& 0& 0& 0& 0& 0& 0& 0& 0& 0& 0 \\
       1& 0& 0& 1& 0& -1& 1& 0& 1& 0& -1& 0& 0& 0& 0& -1& 0& 0& 0& -1& 1& 0 \\
       0& -1& 0& -1& 0& 0& 0& 0& 0& 0& 0& -1& 0& 0& 0& 0& 0& -1& 0& 0& 0& 0
       \\
\end{array}\right)$,

$\left(\begin{array}{cccccccccccccccccccccc}
 0& -1& 0& 0& 0& 0& 0& -1& -1& 0& 1& 0& 0& 0& -1& 0& 1& 0& 1& 1& 0& 1 \\
       -1& 0& 0& -1& 1& 0& -1& 0& 0& 0& 0& 0& -1& 0& -1& 0& 0& 0& -1& 0& 1& 0 \\
       0& 0& 1& 0& 1& 0& -1& 1& 0& 0& 0& 0& 0& 0& -1& 0& 0& 1& -1& 0& 0& 0 \\
       1& -1& 1& 0& 0& -1& 0& 1& 0& 0& -1& -1& 0& -1& 0& -1& 0& 0& 0& -1& 0& 0 \\
       0& -1& 0& 0& 0& 1& 0& 0& -1& 0& 0& 0& 1& 0& 0& 0& -1& 1& 0& 1& 0& 0 \\
       0& 0& 0& 0& 0& 0& 0& 0& 0& 0& 0& 0& 1& 0& 0& 0& 0& 0& 0& 0& 0& 0 \\
       -1& 0& 0& -1& 1& -1& 0& 0& 1& 0& 0& 0& -1& 0& 0& 0& 0& -1& -1& -1& 1& 0 \\
       -1& 0& 0& 0& 0& -1& 0& 0& 0& 0& 1& 0& 0& 0& 0& 0& 1& -1& 1& 0& 0& 1 \\
       0& 0& 0& 0& 0& 0& 0& 0& 1& 0& 0& 0& 0& 0& 0& 0& 0& 0& 0& 0& 0& 0 \\
       0& 0& 0& 0& 0& 0& 0& 0& 0& 1& 0& 0& 0& 0& 0& 0& 0& 0& 0& 0& 0& 0 \\
       0& -1& 1& -1& 1& -1& 0& 0& 0& 0& 0& -1& -1& 0& -1& 0& 1& -1& 0& 0& 0& 0 \\
       0& 0& 0& 1& -1& 1& 0& 0& -1& 0& 1& 1& 1& 0& 0& 0& 0& 1& 1& 1& -1& 1 \\
       0& 0& 0& 0& 0& 1& 0& 0& 0& 0& 0& 0& 0& 0& 0& 0& 0& 0& 0& 0& 0& 0 \\
       0& 1& -1& 1& -1& 0& 1& -1& 0& 0& 1& 1& 1& 1& 1& 0& 0& 0& 1& 0& 0& 0 \\
       1& -1& 1& 0& 0& 0& 0& 1& -1& 0& -1& -1& 1& 0& 0& 0& 0& 1& 0& 1& -1& 0 \\
       0& 0& 0& 0& 0& 0& 0& 0& 0& 0& 0& 0& 0& 0& 0& 1& 0& 0& 0& 0& 0& 0 \\
       1& 0& 1& 0& -1& 0& 0& 1& 0& 0& -1& -1& 1& 0& 1& 0& 0& 0& 0& 0& -1& 0 \\
       0& 1& 0& 0& 0& 0& 0& 0& 1& 0& -1& 0& 0& 1& 1& 1& -1& 0& -1& 0& 0& -1 \\
       -1& 1& -1& 1& 0& 1& -1& 0& 0& 0& 1& 1& 1& 1& 0& 1& -1& 1& 0& 1& 0& 0 \\
       0& -1& 1& -1& 1& -1& -1& 1& 0& 0& 0& -1& -1& -1& -1& -1& 1& 0& 0& 0& 0& 1 \\
       0& 0& 1& 0& 1& -1& -1& 1& 0& 0& 0& -1& 0& 0& -1& 0& 1& 0& 0& 0& 0& 1 \\
       0& 0& -1& 0& 0& 0& 1& -1& 0& 0& 0& 1& 0& 0& 0& 0& 0& 0& 0& 0& 1& 0
       \\
\end{array}\right)$,

$\left(\begin{array}{cccccccccccccccccccccc}
0& -1& 1& -1& 1& -1& 0& 0& 0& 0& 0& -1& -1& 0& -1& 0& 1& -1& 0& 0& 0& 0 \\
       1& -1& 1& 0& 0& 0& 0& 0& -1& -1& 0& -1& 0& -1& -1& -1& 1& 0& 1& 0& -1& 0 \\
       0& 0& -1& 0& 0& 0& 1& 0& 0& 0& 0& 0& 0& 0& 1& -1& 0& 0& 1& 0& 0& 0 \\
       0& 1& -1& 1& -1& 1& 0& 0& 0& 1& 0& 1& 1& 0& 1& 1& -1& 1& 0& 1& -1& 0 \\
       0& 0& 0& 1& 0& -1& 0& 0& 0& 0& 0& 0& 0& 0& 0& -1& 0& 0& 1& 0& 0& 1 \\
       0& 1& 0& 1& 0& 0& -1& 1& 1& 0& -1& 0& 0& 0& 0& 0& -1& 1& -1& 0& 0& 0 \\
       0& 0& -1& 1& 0& 0& 1& -1& -1& -1& 1& 0& 1& 1& 0& -1& 0& 0& 1& 0& 0& 0 \\
       0& 0& 0& 0& 0& 0& 0& 0& 1& 1& -1& 0& 0& 0& 0& 1& 0& 0& -1& 0& 0& 0 \\
       -1& 1& -1& 0& 0& -1& 1& 0& 1& 0& 0& 0& 0& 1& 1& 0& 0& -1& 0& -1& 1& 0 \\
       -1& 0& 0& 0& 1& 0& -1& 1& 0& 0& 0& 0& 0& 0& 0& 0& -1& 1& 0& 1& 0& 0 \\
       0& 0& -1& 0& 0& 0& 1& 0& 0& 0& 0& 0& 1& 1& 1& 0& -1& 0& 0& 0& 0& -1 \\
       -1& 0& 0& -1& 1& -1& 0& 0& 1& 1& 0& 0& -1& 0& 0& 0& 1& -1& 0& 0& 1& 1 \\
       0& 0& 0& 0& 0& -1& 0& 0& 1& 0& -1& 0& -1& -1& 0& -1& 0& 0& 0& -1& 1& 0 \\
       -1& 0& 1& -1& 1& 0& -1& 1& 0& 0& 0& -1& 0& 1& 0& 1& 0& 0& -1& 1& 0& 0 \\
       1& 0& 0& 0& 0& 0& 0& 0& 0& 0& 0& 0& 0& 0& 0& 0& 0& 0& 0& 0& 0& 0 \\
       1& 0& 0& 1& -1& 0& 1& -1& 0& 0& 0& 1& 0& 0& 0& 0& 0& 0& 0& -1& 0& 0 \\
       0& 1& 0& 0& 0& 0& -1& 0& 1& 0& 0& 0& 0& 0& 0& 1& 0& 0& -1& 0& 0& 0 \\
       0& 0& 0& 0& 0& 0& 0& 0& -1& -1& 1& 0& 1& 1& 0& 0& 0& 0& 0& 0& 0& 0 \\
       0& 0& 1& 0& 0& -1& 0& 1& 1& 1& -1& 0& -1& -1& 0& 0& 1& 0& 0& -1& 0& 1 \\
       0& 0& -1& 0& 0& 0& 1& -1& 0& 0& 1& 0& 0& 0& 0& 0& 1& -1& 1& 0& 0& 0 \\
       1& 0& 0& 0& 0& 0& 0& 0& 1& 0& -1& 0& -1& -1& 0& 0& 0& 0& -1& -1& 0& -1 \\
       1& -1& 1& 0& 0& 0& 0& 0& -1& -1& 0& -1& 0& 0& -1& -1& 0& 0& 0& 0& 0& 0
       \\
\end{array}\right)$,

$\left(\begin{array}{cccccccccccccccccccccc}
0& -1& 1& 0& 0& 0& -1& 1& 0& 1& -1& 0& 0& -1& -1& 0& 0& 1& 0& 1& 0& 1 \\
       1& -1& 1& -1& 0& 0& 0& 0& 0& 0& -1& -1& 0& 0& 0& 0& 0& 0& -1& 0& 0& -1 \\
       0& 0& 0& -1& 0& 0& 0& 1& 0& 0& -1& -1& 0& 0& 1& 0& -1& 0& -1& 0& 0& -1 \\
       -1& 1& -1& 0& 0& 0& 0& 0& 0& 0& 1& 1& 0& 0& 0& 0& 0& 0& 0& 0& 0& 0 \\
       -1& 0& 0& -1& 1& 0& -1& 0& 0& 0& 0& 0& -1& 0& -1& 0& 0& 0& -1& 0& 1& 0 \\
       0& 1& 0& 0& 0& 0& 0& 0& 0& 0& 0& 0& 0& 1& 0& 1& 0& 0& -1& 0& 0& -1 \\
       0& 0& -1& 0& 0& 0& 1& -1& 0& 0& 0& 1& 0& 0& 0& 0& 0& 0& 0& 0& 1& 0 \\
       0& 0& -1& 1& -1& 0& 1& -1& 0& 0& 0& 1& 0& 0& 0& 0& 0& 0& 1& 0& 0& 0 \\
       0& 0& 0& 0& 0& 0& 0& 0& 0& 0& 0& 0& 0& 0& 0& 1& 0& 0& 0& 0& 0& 0 \\
       0& 0& 1& -1& 0& -1& 0& 1& 1& 1& -1& -1& -1& -1& 0& 0& 1& -1& -1& -1& 0& 0 \\
       0& 0& 0& 0& 0& -1& 0& 1& 1& 1& -1& 0& 0& -1& 0& 0& 0& 0& 0& 0& 0& 1 \\
       0& -1& 0& 0& 0& 0& 0& 0& -1& -1& 0& -1& 0& 0& 0& -1& 0& 0& 1& 0& 0& 0 \\
       0& 0& 1& -1& 0& -1& 0& 0& 0& 0& 0& -1& -1& 0& 0& 0& 1& -1& 0& -1& 0& 0 \\
       0& 0& 0& 1& 0& 0& 0& 0& 0& 0& 0& 0& 0& 0& 0& 0& 0& 0& 0& 0& 0& 0 \\
       0& 0& 0& 0& 0& 1& -1& 1& -1& 0& 0& 0& 1& 0& 0& 0& -1& 1& 0& 1& -1& 0 \\
       -1& 1& -1& 1& 0& 0& 0& 0& 1& 1& 0& 2& 0& 0& 0& 0& -1& 1& 0& 0& 1& 1 \\
       0& 0& 0& 0& 0& 0& 0& 0& -1& 0& 1& 0& 0& 0& 0& 0& 1& 0& 1& 0& -1& 0 \\
       0& 0& 0& 0& 1& 0& 0& 0& 0& 0& 0& 0& 0& 0& 0& 0& 0& 0& 0& 0& 0& 0 \\
       0& 0& 0& 0& 0& 0& 0& 0& 1& 0& -1& 0& -1& 0& 0& 0& -1& 0& -1& -1& 1& -1 \\
       1& -1& 0& 0& -1& 1& 0& 0& -1& 0& 0& 0& 1& -1& 0& 0& 0& 1& 1& 1& -1& 0 \\
       1& 0& 1& 0& -1& 0& 0& 1& 0& 1& -1& 0& 0& -1& 0& 0& 0& 1& 0& 0& -1& 0 \\
       0& 1& 0& 1& 0& 1& -1& 0& 0& 0& 0& 1& 1& 1& 0& 1& -1& 1& -1& 1& 0& 0
       \\
\end{array}\right)$,

$\left(\begin{array}{cccccccccccccccccccccc}
0& 0& -1& 0& 0& 0& 1& 0& 0& 0& 0& 0& 1& 1& 1& 0& -1& 0& 0& 0& 0& -1 \\
       0& 1& -1& 0& -1& 1& 0& 0& 1& 1& -1& 1& 0& 0& 1& 1& -1& 0& -1& 0& 0& -1 \\
       0& 0& 0& 0& 0& 0& 0& 0& 0& 0& 0& 0& 0& 0& 0& 0& 0& 0& 0& 0& 1& 0 \\
       -1& 0& 0& 0& 1& 0& -1& 0& 0& 0& 1& 1& -1& 0& -1& 0& 0& 0& 0& 0& 1& 1 \\
       -1& 0& 0& -1& 1& 0& -1& 1& 0& 1& 0& 0& 0& 0& 0& 1& 0& 0& -1& 1& 0& 0 \\
       0& 0& 0& -1& 0& 0& 0& 0& 0& 0& 0& 0& -1& -1& 0& 0& 1& -1& 0& -1& 0& 0 \\
       0& 0& 1& -1& 1& 0& -1& 1& 0& 0& -1& -1& -1& 0& 0& 0& 0& 0& -1& 0& 0& -1 \\
       0& 0& 0& 0& 1& -1& 0& 0& 0& -1& 0& 0& -1& 0& -1& -1& 0& 0& 0& -1& 1& 0 \\
       1& 0& 1& 0& 0& 0& 0& 0& 0& -1& -1& -1& 0& 0& 0& 0& 0& 0& -1& -1& 0& -1 \\
       0& 0& -1& 1& -1& 0& 1& 0& 0& 1& 0& 1& 1& 0& 1& 0& 0& 0& 1& 0& 0& 1 \\
       0& 1& 0& 1& 0& 0& 0& 1& 0& 0& 0& 0& 1& 1& 1& 0& -1& 1& 0& 0& 0& 0 \\
       0& -1& 0& 0& 1& -1& 0& 0& 0& -1& 0& -1& 0& 0& -1& -1& 0& 0& 0& 0& 1& 0 \\
       0& 0& 0& 0& -1& 1& 0& 0& -1& 0& 0& 0& 1& 0& 0& 1& 0& 0& 0& 1& -1& 0 \\
       1& -1& 0& 0& 0& -1& 1& -1& 0& -1& 0& -1& 0& 0& 0& -1& 1& -1& 1& -1& 0& 0 \\
       0& 0& 0& 0& 0& 0& 0& 0& 0& 0& 1& 0& 0& 0& 0& 0& 0& 0& 0& 0& 0& 0 \\
       -1& 0& 0& -1& 1& 0& 0& 0& -1& 0& 1& 0& 0& 1& 0& 0& 0& 0& 0& 1& 0& 0 \\
       0& 0& 0& 0& 0& 0& 0& -1& 0& -1& 1& 0& -1& 0& -1& 0& 1& -1& 0& -1& 0& 0 \\
       1& 0& 1& 0& 0& 0& 0& 0& 0& 0& 0& -1& 0& 0& 0& 0& 1& 0& 0& 0& -1& 0 \\
       0& 0& -1& 0& 0& 0& 1& -1& 0& 0& 0& 1& 0& 0& 0& 0& 0& 0& 0& 0& 1& 0 \\
       0& 0& 0& 0& 0& 0& 0& 0& 1& 0& -1& 0& -1& 0& 0& 0& -1& 0& -1& -1& 1& -1 \\
       0& 0& -1& 0& -1& 1& 1& -1& -1& 0& 1& 1& 0& 0& 0& 0& 0& 0& 1& 0& 0& 0 \\
       0& 0& 0& -1& 0& 0& 0& 0& 0& 0& 0& -1& 0& 0& 1& 0& 0& -1& 0& 0& -1& -1
       \\
\end{array}%
\right)$.

The character table of $G^{s_2}$:\\
\begin{tabular}{c|cccccccccccccccccccc}
  & & & & & & & & & &10 & & & & & & & & & &20\\\hline
\\$\chi_{2}^{(1)}$&1&1&1&1&1&1&1&1&1&1&1&1&1&1&1&1&1&1&1&1
\\$\chi_{2}^{(2)}$&1&1&-1&1&-1&1&-1&-1&-1&-1&1&1&-1&-1&1&-1&-1&1&1&-1
\\$\chi_{2}^{(3)}$&1&1&C&1&/C&C&1&1&C&1&/C&/C&C&/C&C&C&1&1&C&/C
\\$\chi_{2}^{(4)}$&1&1&-C&1&-/C&C&-1&-1&-C&-1&/C&/C&-C&-/C&C&-C&-1&1&C&-/C
\\$\chi_{2}^{(5)}$&1&1&/C&1&C&/C&1&1&/C&1&C&C&/C&C&/C&/C&1&1&/C&C
\\$\chi_{2}^{(6)}$&1&1&-/C&1&-C&/C&-1&-1&-/C&-1&C&C&-/C&-C&/C&-/C&-1&1&/C&-C
\\
\end{tabular}%

 \noindent  \noindent where
A = E(5); B = E(5)$^2$; C = E(3)
  = (-1+ER(-3))/2 = b3;
D = E(15)$^2$; E = E(15)$^7$; F = E(15)$^{14}$; G = E(15)$^4$.

The generators of $G^{s_3}$ are:\\
 $\left(%
\right)$.

The representatives of conjugacy classes of   $G^{s_3}$ are:\\
$\left(%
\right)$.

The character table of $G^{s_{3}}$:\\

\begin{tabular}{c|cccccccccccccccccccc}
  & & & & & & & & & &10 & & & & & & & & & &20\\\hline
\\$\chi_{3}^{(1)}$&1&1&1&1&1&1&1&1&1&1&1&1&1&1&1&1&1&1&1&1
\\$\chi_{3}^{(2)}$&1&1&-1&1&-1&1&-1&-1&-1&-1&1&1&-1&-1&1&-1&-1&1&1&-1
\\$\chi_{3}^{(3)}$&1&1&C&1&/C&C&1&1&C&1&/C&/C&C&/C&C&C&1&1&C&/C
\\$\chi_{3}^{(4)}$&1&1&-C&1&-/C&C&-1&-1&-C&-1&/C&/C&-C&-/C&C&-C&-1&1&C&-/C
\\$\chi_{3}^{(5)}$&1&1&/C&1&C&/C&1&1&/C&1&C&C&/C&C&/C&/C&1&1&/C&C
\\$\chi_{3}^{(6)}$&1&1&-/C&1&-C&/C&-1&-1&-/C&-1&C&C&-/C&-C&/C&-/C&-1&1&/C&-C
\\
\end{tabular}%

 \noindent  \noindent where
A = E(5)$^3$; B = E(5); C = E(3)
  = (-1+ER(-3))/2 = b3;
D = E(15)$^{11}$; E = E(15); F = E(15)$^2$; G = E(15)$^7$.

The generators of $G^{s_{4}}$ are:\\
$\left(%
\right)$.

The representatives of conjugacy classes of   $G^{s_4}$ are:\\

$\left(%
\right)$.

The character table of $G^{s_{4}}$:\\

\begin{tabular}{c|cccccccccccccccccccc}
  & & & & & & & & & &10 & & & & & & & & & &20\\\hline

\\$\chi_{4}^{(1)}$&1&1&1&1&1&1&1&1&1&1&1&1&1&1&1&1&1&1&1&1
\\$\chi_{4}^{(2)}$&1&1&-1&1&-1&1&-1&-1&-1&-1&1&1&-1&-1&1&-1&-1&1&1&-1
\\$\chi_{4}^{(3)}$&1&1&C&1&/C&-C&-1&-1&C&-1&-/C&-/C&C&/C&-C&C&-1&1&-C&/C
\\$\chi_{4}^{(4)}$&1&1&/C&1&C&-/C&-1&-1&/C&-1&-C&-C&/C&C&-/C&/C&-1&1&-/C&C
\\$\chi_{4}^{(5)}$&1&1&-/C&1&-C&-/C&1&1&-/C&1&-C&-C&-/C&-C&-/C&-/C&1&1&-/C&-C
\\$\chi_{4}^{(6)}$&1&1&-C&1&-/C&-C&1&1&-C&1&-/C&-/C&-C&-/C&-C&-C&1&1&-C&-/C
\\
\end{tabular}%

 \noindent  \noindent where
A = E(5)$^3$; B = E(5); C = -E(3)$^2$
  = (1+ER(-3))/2 = 1+b3;
D = -E(15)$^2$; E = -E(15)$^4$; F = -E(15)$^8$; G = -E(15)$^14$.

The generators of $G^{s_{5}}$ are:\\

$\left(%
\right)$.

The representatives of conjugacy classes of   $G^{s_{5}}$ are:\\

$\left(%
\right)$.

The character table of $G^{s_{5}}$:\\

%

 \noindent  \noindent where
A = E(5)$^3$; B = E(5); C = -E(3)$^2$
  = (1+ER(-3))/2 = 1+b3;
D = -E(15)$^2$; E = -E(15)$^4$; F = -E(15)$^8$; G = -E(15)$^{14}$.

The generators of $G^{s_{6}}$ are:\\
$\left(%
\right)$.

The representatives of conjugacy classes of   $G^{s_{6}}$ are:\\
$\left(%
\right)$.

The character table of $G^{s_{6}}$:\\

%

 \noindent
 \noindent where
A = 5*E(5)$^4$; B = 5*E(5)$^3$; C = -2*E(5)-2*E(5)$^4$
  = 1-ER(5) = 1-r5;
D = E(5)+2*E(5)$^2$+2*E(5)$^3$+E(5)$^4$
  = (-3-ER(5))/2 = -2-b5;
E = E(3)$^2$
  = (-1-ER(-3))/2 = -1-b3;
F = -E(5)$^4$; G = -E(5)$^3$; H = -E(15)$^7$; I = -E(15)$^2$; J =
-E(15)$^4$; K = -E(15)$^{14}$.

The generators of $G^{s_{7}}$ are:\\
$\left(%
\right)$.

The representatives of conjugacy classes of   $G^{s_{7}}$are:\\
$\left(%
\right)$.

The character table of $G^{s_{7}}$:\\

%

 \noindent  \noindent where
A = 9*E(3)$^2$
  = (-9-9*ER(-3))/2 = -9-9b3;
B = 18*E(3)$^2$
  = -9-9*ER(-3) = -9-9i3;
C = 27*E(3)$^2$
  = (-27-27*ER(-3))/2 = -27-27b3;
D = 36*E(3)$^2$
  = -18-18*ER(-3) = -18-18i3;
E = 45*E(3)$^2$
  = (-45-45*ER(-3))/2 = -45-45b3;
F = 54*E(3)$^2$
  = -27-27*ER(-3) = -27-27i3;
G = E(3)$^2$
  = (-1-ER(-3))/2 = -1-b3;
H = -2*E(3)$^2$
  = 1+ER(-3) = 1+i3;
I = 3*E(3)$^2$
  = (-3-3*ER(-3))/2 = -3-3b3;
J = 4*E(3)$^2$
  = -2-2*ER(-3) = -2-2i3;
K = 5*E(3)$^2$
  = (-5-5*ER(-3))/2 = -5-5b3;
L = -6*E(3)$^2$
  = 3+3*ER(-3) = 3+3i3;
M = -E(3)+2*E(3)$^2$
  = (-1-3*ER(-3))/2 = -2-3b3;
N = -E(5)-E(5)$^4$
  = (1-ER(5))/2 = -b5;
O = -E(15)-E(15)$^4$; P = -E(15)$^7$-E(15)$^{13}$.

The generators of $G^{s_{8}}$ are:\\
$\left(
\right)$.

The representatives of conjugacy classes of   $G^{s_{8}}$
are:\\
$\left(%
\right)$.

The character table of $G^{s_{8}}$:\\

%

 \noindent  \noindent where
A = E(5)$^3$; B = E(5); C = -E(3)$^2$
  = (1+ER(-3))/2 = 1+b3;
D = -E(15)$^2$; E = -E(15)$^4$; F = -E(15)$^8$; G = -E(15)$^{14}$.

The generators of $G^{s_{9}}$ are:\\

$\left(%
\right)$.

The representatives of conjugacy classes of   $G^{s_{9}}$
are:\\
$\left(%
\right)$.

The character table of $G^{s_{9}}$:\\

%

 \noindent  \noindent where
A = -E(3)+E(3)$^2$
  = -ER(-3) = -i3;
B = -E(7)-E(7)$^2$-E(7)$^4$
  = (1-ER(-7))/2 = -b7;
C = -E(15)-E(15)$^2$-E(15)$^4$-E(15)$^8$
  = (-1-ER(-15))/2 = -1-b15.

The generators of $G^{s_{10}}$ are:\\
$\left(%
\right)$.

The representatives of conjugacy classes of   $G^{s_{10}}$
are:\\
$\left(%
\right)$.

  The character table of $G^{s_{{10}}}$:\\
\setlength{\tabcolsep}{4pt}
%

 \noindent  \noindent where
A = E(3)
  = (-1+ER(-3))/2 = b3;
B = -E(5)-E(5)$^4$
  = (1-ER(5))/2 = -b5;
C = -E(15)$^2$-E(15)$^8$; D = -E(15)$^{11}$-E(15)$^{14}$; E = 2*E(3)
  = -1+ER(-3) = 2b3;
F = 3*E(3)
  = (-3+3*ER(-3))/2 = 3b3;
G = 4*E(3)
  = -2+2*ER(-3) = 4b3;
H = 5*E(3)
  = (-5+5*ER(-3))/2 = 5b3;
I = 6*E(3)
  = -3+3*ER(-3) = 6b3.

The generators of $G^{s_{11}}$ are:\\
$\left(%
\right)$.

 The character table of $G^{s_{{11}}}$:\\

%

 \noindent  \noindent where
A = E(7)$^5$; B = E(7)$^3$; C = E(7).

The generators of $G^{s_{12}}$ are:\\
$\left(%
\right)$,

The representatives of conjugacy classes of   $G^{s_{12}}$
are:\\
$\left(%
\right)$.

The character table of $G^{s_{{12}}}$:\\

%

 \noindent  \noindent where

A = E(7)$^4$; B = E(7); C = E(7)$^5$.

The generators of $G^{s_{13}}$ are:\\
$\left(%
\right)$.

The character table of $G^{s_{{13}}}$:\\

%

 \noindent  \noindent where
A = E(7)$^2$; B = E(7)$^4$; C = E(7)$^6$.

The generators of $G^{s_{14}}$ are:\\
$\left(%
\right)$.

The character table of $G^{s_{{14}}}$:\\
%

 \noindent  \noindent where
A = E(7)$^2$; B = E(7)$^4$; C = E(7)$^6$.

The generators of $G^{s_{15}}$ are:\\
$\left(%
\right)$.

The character table of $G^{s_{{15}}}$:\\

%

 \noindent  \noindent where
A = E(3)$^2$
  = (-1-ER(-3))/2 = -1-b3;
B = -E(4)
  = -ER(-1) = -i;
C = -E(12)$^{11}$.

The generators of $G^{s_{16}}$ are:\\
$\left(%
\right)$.

The character table of $G^{s_{{16}}}$:\\

%

 \noindent  \noindent where
A = -E(4)
  = -ER(-1) = -i;
B = -1+E(4)
  = -1+ER(-1) = -1+i;
C = -2*E(4)
  = -2*ER(-1) = -2i;
D = -4*E(4)
  = -4*ER(-1) = -4i.

The generators of $G^{s_{17}}$ are:\\
$\left(%
\right)$.

The character table of $G^{s_{{17}}}$:\\
\setlength{\tabcolsep}{4pt}
%
\noindent  \noindent where

A = E(3)$^2$
  = (-1-ER(-3))/2 = -1-b3;
B = 3*E(3)$^2$
  = (-3-3*ER(-3))/2 = -3-3b3;
C = 4*E(3)$^2$
  = -2-2*ER(-3) = -2-2i3;
D = 6*E(3)$^2$
  = -3-3*ER(-3) = -3-3i3;
E = 12*E(3)$^2$
  = -6-6*ER(-3) = -6-6i3;
F = -2*E(3)$^2$
  = 1+ER(-3) = 1+i3;
G = -E(3)+2*E(3)$^2$
  = (-1-3*ER(-3))/2 = -2-3b3;
H = -3*E(3)-2*E(3)$^2$
  = (5-ER(-3))/2 = 2-b3;
I = 3*E(3)+E(3)$^2$
  = -2+ER(-3) = -2+i3.

The generators of $G^{s_{18}}$ are:\\
$\left(%
\right)$.

The character table of $G^{s_{{18}}}$:\\
\setlength{\tabcolsep}{4pt}
%

\noindent  \noindent where A = E(3)
  = (-1+ER(-3))/2 = b3
B = E(9)$^7$; C = E(9)$^5$; D = -E(9)$^4$-E(9)$^7$.

The generators of $G^{s_{19}}$ are:\\
$\left(%
\right)$.

The character table of $G^{s_{{19}}}$:\\
\setlength{\tabcolsep}{4pt}
%

\noindent  \noindent where A = E(3)
  = (-1+ER(-3))/2 = b3;
B = E(9)$^7$; C = E(9)$^5$; D = -E(9)$^4$-E(9)$^7$.

The generators of $G^{s_{20}}$ are:\\
$\left(%
\right)$.

The character table of $G^{s_{{20}}}$:\\
%

 \noindent  \noindent where
A = E(11)$^2$; B = E(11)$^4$; C = E(11)$^6$; D = E(11)$^8$; E =
E(11)$^{10}$.

The generators of $G^{s_{21}}$ are:\\
$\left(%
\right)$.

The character table of $G^{s_{{21}}}$:\\
%

 \noindent  \noindent where
A = E(11); B = E(11)$^2$; C = E(11)$^3$; D = E(11)$^4$; E =
E(11)$^5$.

The generators of $G^{s_{22}}$ are:\\
$\left(%
\right)$.

The character table of $G^{s_{{22}}}$:\\
%

 \noindent  \noindent where
A = E(4)
  = ER(-1) = i;
B = E(8).

The generators of $G^{s_{23}}$ are:\\
$\left(%
\right)$.

The character table of $G^{s_{23}}$:\\
%

 \noindent  \noindent where
A = -E(3)
  = (1-ER(-3))/2 = -b3;
B = 2*E(3)$^2$
  = -1-ER(-3) = -1-i3.

The generators of $G^{s_{24}}$ are:\\
$\left(%
\right)$.

The character table of $G^{s_{24}}$:\\
\setlength{\tabcolsep}{4pt}
%

\noindent where A = E(5)$^3$; B = E(5).

    \vskip 0.3cm

  {\large\bf Acknowledgement}: We would like to thank Prof.
N. Andruskiewitsch and Dr. F. Fantino for suggestions and help.

\end {document}